\documentclass[a4paper]{article}
\usepackage{geometry}
\usepackage{amssymb}
\usepackage{latexsym}
\usepackage{amsmath}
\usepackage{amsthm} 
\usepackage{indentfirst}
\usepackage{graphicx}
\usepackage[colorlinks=true]{hyperref}
\usepackage{mathrsfs}
\usepackage{chngcntr}
\usepackage{color}

\usepackage{tikz}

\newtheorem{theorem}{Theorem}[section]

\newtheorem{lemma}{Lemma}[section]
\newtheorem{Pro}{Proposition}[section]
\newtheorem{Exam}{Example}[section]
\newtheorem{remark}{Remark}[section]
\newtheorem{Def}{Definition}[section]

\numberwithin{equation}{section}

\begin{document}

\title{Bound states of nonlinear Dirac equations on periodic quantum graphs}

\author{
  Zhipeng Yang\thanks{Email: yangzhipeng326@163.com}
  \quad and\quad
  Ling Zhu\\[2pt]
  \small Yunnan Key Laboratory of Modern Analytical Mathematics and Applications, Kunming, China\\
  \small Department of Mathematics, Yunnan Normal University, Kunming, China\\
}

\date{}
\maketitle

\begin{abstract}
We study nonlinear Dirac equations (NLDE) on periodic quantum graphs endowed with Kirchhoff-type vertex conditions. Our main goal is to establish existence and multiplicity of bound states, which arise as critical points of the associated NLDE action functional. The underlying Dirac operator has a spectral gap around the origin, so the corresponding functional is strongly indefinite, and in addition the Palais--Smale condition fails due to the noncompactness and the periodic structure of the graph. To overcome these difficulties, we combine the spectral properties of the periodic Dirac operator with critical point theorems for strongly indefinite functionals and a concentration--compactness analysis adapted to periodic quantum graphs, and derive existence and multiplicity results for bound states with frequencies lying in the spectral gap.
\end{abstract}




	\section{Introduction and main results}
	
	The study of quantum graphs has emerged as a vibrant field at the intersection of
	mathematical physics, spectral theory, and nonlinear analysis \cite{MR3013208}.
	Quantum graphs---metric graphs equipped with differential operators (Hamiltonians)
	and vertex conditions---serve as idealized models for complex systems such as
	photonic crystals, carbon nanostructures, and quantum networks
	\cite{MR3362506,MR2369953}. A paradigmatic example of a Hamiltonian is the
	Laplacian $-\Delta_{\mathcal G}$ acting along the edges with Kirchhoff conditions at the vertices (see Section~2 for a precise
	definition). These structures inherit both the continuous nature of differential
	equations and the discrete combinatorial features of graphs, offering a rich
	framework to explore phenomena such as wave propagation, spectral gaps, and
	localized modes. A central challenge in this domain lies in understanding how the
	geometry of the graph and the vertex conditions influence the existence and
	properties of solutions to nonlinear evolution equations \cite{MR2459876}.
	
	In this paper we are interested in nonlinear Dirac equations on periodic quantum
	graphs. More precisely, we assume that $\mathcal G$ is a connected metric graph
	carrying a free, cocompact action of the group $\mathbb Z^{d}$ by graph
	automorphisms. Equivalently, there exists a compact connected subgraph
	$\mathcal K\subset\mathcal G$ (a fundamental cell) such that
	\[
	\mathcal G = \bigcup_{k\in\mathbb Z^{d}} T^{k}(\mathcal K),
	\]
	where each $T^{k}$ is a graph isomorphism and the translates
	$T^{k}(\mathcal K)$ intersect only along boundary vertices whenever
	$k\neq\ell$; see Section~2 for details. On such graphs, the Dirac operator acts
	as a first-order system along the edges and is coupled by Kirchhoff-type
	conditions at the vertices.
	
	On the other hand, the Dirac equation, originally formulated to describe
	relativistic electrons, has gained renewed interest in condensed matter physics
	due to its relevance for materials with linear dispersion relations, such as
	graphene \cite{castro2009electronic,katsnelson2007graphene,thaller1992dirac}.
	While the linear Dirac equation has been extensively studied, its nonlinear
	counterpart introduces self-interaction terms that model rich phenomena such as
	solitons, localized particle-like states, and nonlinear optical effects in
	relativistic quantum systems \cite{MR1700174}. On quantum graphs, the Dirac
	operator acts as a first-order system on each edge, coupling the spinor
	components and interacting with the underlying graph topology
	\cite{MR3557323,MR1965289,MR2459876}.
	
	In particular, in the simplified setting of the infinite 3-star graph, the authors of \cite{MR3871194} proposed the study of
	nonlinear Dirac equations on networks, where the Dirac operator on each edge is
	given by
	\begin{equation}\label{eq-1.1}
		\mathcal D
		= -ic\,\frac{d}{dx}\otimes\sigma_{1}
		+ mc^{2}\otimes\sigma_{3},
	\end{equation}
	where $m>0$ represents the mass of the particle and $c>0$ is the speed of light.
	Here $\sigma_{1}$ and $\sigma_{3}$ are the Pauli matrices
	\[
	\sigma_{1} = \begin{pmatrix} 0 & 1\\ 1 & 0\end{pmatrix},
	\qquad
	\sigma_{3} = \begin{pmatrix} 1 & 0\\ 0 & -1\end{pmatrix}.
	\]
	As is well known, solutions to the Dirac equation \eqref{eq-1.1} are spinors
	$\chi = (\chi^{1},\chi^{2})^{T}$. Motivated by pure power self-interactions, it is natural to look for stationary solutions of a nonlinear Dirac equation of the form
\[
\chi(t,x) = e^{i\omega t} u(x),\qquad \omega\in\mathbb R,
\]
which leads to the stationary nonlinear Dirac equation
\begin{equation}\label{eq-1.2}
\mathcal D u + \omega u = |u|^{p-2}u.
\end{equation}
	
	Subsequently, the authors of \cite{MR3934110} initiated a systematic study of
	bound states and the nonrelativistic limit for NLDE on noncompact quantum graphs.
	To treat more intricate graph topologies, they considered Kirchhoff-type
	extensions of the Dirac operator on general metric graphs. A crucial feature of
	their model is a localized nonlinearity, which leads to the equation
	\begin{equation}\label{eq-1.3}
		\mathcal D\psi+\omega\psi=\chi_{\mathcal K}|\psi|^{p-2}\psi,
	\end{equation}
	where $\chi_{\mathcal K}$ denotes the characteristic function of a compact
	subgraph $\mathcal K\subset\mathcal G$. Furthermore, they rigorously proved that,
	in the nonrelativistic limit $c\to+\infty$, the bound states of
	\eqref{eq-1.3} converge to those of the associated nonlinear Schr\"odinger
	equation on the same graph,
	\[
	-\Delta_{\mathcal G}u+\omega u=\chi_{\mathcal K}|u|^{p-2}u.
	\]
	
	It is worth noting that in \eqref{eq-1.3}, as well as in the equation analyzed in
	\cite{MR4200758}, the nonlinearity is a pure power term, which is manifestly
	non-covariant. Such nonlinearities are, however, quite natural in nonlinear
	optics. From a theoretical standpoint this does not create a conceptual
	difficulty, since the nonlinear Dirac equation is understood as an effective
	model rather than a fully covariant field theory. Motivated by this point of
	view, in the present work we consider more general nonlinearities, extending the
	analysis beyond the pure power case.
	
	More precisely, we study the stationary NLDE
	\begin{equation}\label{eq-1.4}
		\mathcal D u + \omega u
		= F_{u}(x,u)
		\quad\text{on }\mathcal G,
	\end{equation}
	where $F_{u}$ denotes the gradient of $F$ with respect to $u\in\mathbb C^{2}$.
	We assume that the underlying graph $\mathcal G$ is a noncompact periodic quantum
	graph as described above, and that the nonlinearity $F$ is $\mathbb Z^{d}$–periodic
	in the spatial variable.
	
	For the nonlinear term we impose the following conditions:
	\begin{itemize}
		\item[$(F_{0})$] $F\in C^{1}(\mathcal G\times\mathbb C^{2},[0,+\infty))$.
		\item[$(F_{1})$] $F(x,u)$ is $\mathbb Z^{d}$-periodic in the spatial variable
		in the sense that
		\[
		F(T^{k}x,u) = F(x,u)
		\quad\text{for all }k\in\mathbb Z^{d},\ x\in\mathcal G,\ u\in\mathbb C^{2},
		\]
		where $T^{k}$ denotes the graph translations.
		\item[$(F_{2})$] $F_{u}(x,u) = o(|u|)$ as $u\to 0$, uniformly in $x\in\mathcal G$.
	\end{itemize}
	We also set
	\[
	\omega_{0}
	= \min\{mc^{2}+\omega,\; mc^{2}-\omega\},
	\qquad
	\hat F(x,u)
	= \frac{1}{2}F_{u}(x,u)\cdot u - F(x,u).
	\]
	Furthermore we require:
	\begin{itemize}
		\item[$(F_{3})$] There exists $b>mc^{2}+\omega$ such that
\[
\frac{|F_{u}(x,u)-bu|}{|u|}\to 0\quad\text{as }|u|\to\infty,
\]
uniformly in $x\in\mathcal G$.
Moreover, this same $b$ satisfies $b-\omega\notin\sigma_{p}(\mathcal D)$.

		\item[$(F_{4})$] $\hat F(x,u)\ge 0$ for all $(x,u)$, and there exists
		$\delta_{1}\in(0,\omega_{0})$ such that
		\[
		\hat F(x,u)\ge\delta_{1}
		\quad\text{whenever}\quad
		|F_{u}(x,u)|\ge (\omega_{0}-\delta_{1})|u|.
		\]
		\item[$(F_{5})$] $F\in C^{2}(\mathcal G\times\mathbb C^{2},[0,+\infty))$ and
		there exist $\nu\in[0,1)$ and $c_{1}>0$ such that
		\[
		|F_{uu}(x,u)|\le c_{1}\bigl(1+|u|^{\nu}\bigr)
		\quad\text{for all }(x,u)\in\mathcal G\times\mathbb C^{2}.
		\]
	\end{itemize}
	
	Our main result can now be stated as follows.
	
	\begin{theorem}\label{theo-1.1}
		Let $\mathcal G$ be a noncompact periodic quantum graph carrying a free,
		cocompact action of $\mathbb Z^{d}$. Assume that $m,c>0$, $\omega\in(-mc^{2},mc^{2})$
		and that $(F_{0})$--$(F_{5})$ are satisfied. Then the nonlinear Dirac equation
		\eqref{eq-1.4} admits at least one bound state $u$.
		
		In addition, if $F(x,u)$ is even in $u$, then
		\eqref{eq-1.4} admits infinitely many geometrically distinct bound states.
	\end{theorem}
	
	\begin{remark}\label{rem-1.1}
		Two bound states $u_{1}$ and $u_{2}$ are said to be geometrically distinct if
		\[
		u_{2}\neq k*u_{1}
		\quad\text{for all }k\in\mathbb Z^{d},
		\]
		where the action of $\mathbb Z^{d}$ is given by
		\[
		(k*u)(x):=u\bigl(T^{k}x\bigr),
		\qquad x\in\mathcal G.
		\]
	\end{remark}
	
	A distinctive feature of the NLDE is the strong indefiniteness of its action
	functional. Unlike the Schr\"odinger case, where the associated functional is typically bounded from below \cite{MR1400007}, 
the Dirac functional is never coercive in any natural Hilbert space, because the spectrum of the Dirac operator is unbounded both above and below (see Section~2).
This indefiniteness precludes direct minimization techniques and calls for more
	sophisticated critical point theories \cite{MR2389415}. Moreover, the
	noncompactness and periodicity of the graph undermine the global
	Palais--Smale condition, a cornerstone of classical variational methods.
	Consequently, standard tools such as the mountain pass theorem or symmetric
	minimax principles must be substantially adapted or replaced.
	
	In this work we employ a combination of these strategies. First, we decompose the
	underlying Hilbert space into the positive and negative spectral subspaces of
	the periodic Dirac operator, thereby isolating the contributions of the strongly
	indefinite linear part. Next, we exploit the periodicity of both the graph and
	the nonlinearity together with the Kirchhoff-type vertex conditions to construct
	critical points of the action functional by topological and variational
	arguments. Crucially, we do not assume a global Palais--Smale condition; instead
	we verify suitable compactness properties along carefully chosen sequences that
	reflect the geometry of the periodic graph and the translation invariance of the
	problem.
	
	The paper is organized as follows. In Section~2 we recall the definition of
	periodic quantum graphs, introduce the Dirac operator with Kirchhoff-type vertex
	conditions, and formulate the stationary NLDE as a variational problem. Section~3
	is devoted to the proofs of the existence and multiplicity results stated in
	Theorem~\ref{theo-1.1}.

\section{Preliminaries}

\subsection{Periodic quantum graphs and functional setting}

We briefly recall the basic definitions and functional setting for quantum
graphs, referring to \cite{MR3385179,MR3013208} and the references therein
for further details. Throughout the paper, integrals on $\mathcal G$ are taken
with respect to the one-dimensional Lebesgue measure along the edges.

A (metric) graph is a pair $\mathcal G=(\mathcal V,\mathcal E)$, where
$\mathcal V$ is the set of vertices and $\mathcal E$ is the set of edges.
To each edge $e\in\mathcal E$ we associate either a bounded interval
$I_{e}=[0,\ell_{e}]$ with length $\ell_{e}>0$ or a half-line
$I_{e}=[0,+\infty)$, together with an orientation and a coordinate
$x_{e}\in I_{e}$. The endpoints of the intervals are identified with
vertices in $\mathcal V$ according to the combinatorial structure of the
graph, so that $\mathcal G$ becomes a one-dimensional metric space obtained
by gluing the intervals $I_{e}$ at their endpoints.

In this work we consider quantum graphs that are periodic under a free,
cocompact action of the lattice $\mathbb Z^{d}$ by graph automorphisms.
More precisely, we assume that $\mathcal G$ is a connected metric graph such that
there exists a group homomorphism
\[
\mathbb Z^{d}\ni k \longmapsto T^{k}\in\operatorname{Aut}(\mathcal G)
\]
with the following properties:
\begin{itemize}
	\item[$(a)$] the action $\mathbb Z^{d}\times\mathcal G\to\mathcal G$,
	$(k,x)\mapsto T^{k}(x)$, is free and by isometries on each edge;
	\item[$(b)$] the quotient graph $\mathcal G/\mathbb Z^{d}$ is compact
	(a finite metric graph);
	\item[$(c)$] there exists a compact connected subgraph
	$\mathcal K\subset\mathcal G$ such that
	\[
	\mathcal G = \bigcup_{k\in\mathbb Z^{d}} T^{k}(\mathcal K),
	\]
	and $T^{k}(\mathcal K)\cap T^{\ell}(\mathcal K)$ consists only of boundary
	vertices whenever $k\neq\ell$.
\end{itemize}
The set $\mathcal K$ is called a fundamental cell of the periodic graph
$\mathcal G$. Notice that each edge of $\mathcal G$ is contained in some
translate $T^{k}(\mathcal K)$, and that $\mathcal G$ has no half-lines: all
edges are finite line segments, but there are infinitely many of them.

A function $u:\mathcal G\to\mathbb C$ can be identified with a family
$(u_{e})_{e\in\mathcal E}$, where $u_{e}:I_{e}\to\mathbb C$ is the
restriction of $u$ to the edge $I_{e}$. For $1\le p<\infty$ we consider
\[
L^{p}(\mathcal G)
=\Bigl\{u=(u_{e})_{e\in\mathcal E}:\ u_{e}\in L^{p}(I_{e})\ \text{for all }e,\
\sum_{e\in\mathcal E}\|u_{e}\|_{L^{p}(I_{e})}^{p}<\infty\Bigr\},
\]
with norm
\[
\|u\|_{L^{p}(\mathcal G)}^{p}
=\sum_{e\in\mathcal E}\|u_{e}\|_{L^{p}(I_{e})}^{p},
\qquad 1\le p<\infty,
\]
and
\[
\|u\|_{L^{\infty}(\mathcal G)}
=\sup_{e\in\mathcal E}\|u_{e}\|_{L^{\infty}(I_{e})}.
\]
We define
\[
H^{1}(\mathcal G)
=\Bigl\{u=(u_{e})_{e\in\mathcal E}:\ u_{e}\in H^{1}(I_{e})\ \text{for all }e,\
\sum_{e\in\mathcal E}\|u_{e}\|_{H^{1}(I_{e})}^{2}<\infty\Bigr\},
\]
with norm
\[
\|u\|_{H^{1}(\mathcal G)}^{2}
=\|u'\|_{L^{2}(\mathcal G)}^{2} + \|u\|_{L^{2}(\mathcal G)}^{2},
\]
where $u'=(u_{e}')_{e\in\mathcal E}$ denotes the family of weak derivatives
along each edge.

A spinor $u:\mathcal G\to\mathbb C^{2}$ is a map $u=(u^{1},u^{2})^{T}$ whose
components $u^{1},u^{2}:\mathcal G\to\mathbb C$ are scalar functions. Equivalently,
one may regard $u$ as a family of $2$-spinors
\[
u_{e} = \binom{u_{e}^{1}}{u_{e}^{2}} : I_{e}\to\mathbb C^{2},
\qquad e\in\mathcal E.
\]
We consider
\[
L^{p}(\mathcal G,\mathbb C^{2})
=\{u=(u^{1},u^{2})^{T}:\ u^{1},u^{2}\in L^{p}(\mathcal G)\},
\]
with
\[
\|u\|_{L^{p}(\mathcal G,\mathbb C^{2})}^{p}
=\|u^{1}\|_{L^{p}(\mathcal G)}^{p}
+\|u^{2}\|_{L^{p}(\mathcal G)}^{p},
\quad 1\le p<\infty,
\]
and
\[
\|u\|_{L^{\infty}(\mathcal G,\mathbb C^{2})}
=\max\bigl\{\|u^{1}\|_{L^{\infty}(\mathcal G)},
\|u^{2}\|_{L^{\infty}(\mathcal G)}\bigr\}.
\]
Similarly,
\[
H^{1}(\mathcal G,\mathbb C^{2})
=\{u=(u^{1},u^{2})^{T}:\ u^{1},u^{2}\in H^{1}(\mathcal G)\},
\]
with
\[
\|u\|_{H^{1}(\mathcal G,\mathbb C^{2})}^{2}
=\|u^{1}\|_{H^{1}(\mathcal G)}^{2}
+\|u^{2}\|_{H^{1}(\mathcal G)}^{2}.
\]

On each finite edge $I_{e}$ the one-dimensional Sobolev embedding
$H^{1}(I_{e})\hookrightarrow C(I_{e})\cap L^{\infty}(I_{e})$ yields
\[
\|u_{e}\|_{L^{\infty}(I_{e})}\le C(\ell_{e})\,\|u_{e}\|_{H^{1}(I_{e})},
\qquad u_{e}\in H^{1}(I_{e}),
\]
where $C(\ell_{e})$ depends only on the length $\ell_{e}$.
Since $\mathcal G/\mathbb Z^{d}$ is a finite metric graph, there are only finitely
many edge lengths in $\mathcal G$, hence $\sup_{e\in\mathcal E}C(\ell_{e})<\infty$.
Moreover, for every $e\in\mathcal E$ we have
$\|u_{e}\|_{H^{1}(I_{e})}\le \|u\|_{H^{1}(\mathcal G)}$ because
$\|u\|_{H^{1}(\mathcal G)}^{2}=\sum_{e\in\mathcal E}\|u_{e}\|_{H^{1}(I_{e})}^{2}$.
Therefore there exists $C_{\infty}>0$ such that
\[
\|u\|_{L^{\infty}(\mathcal G)}
=\sup_{e\in\mathcal E}\|u_{e}\|_{L^{\infty}(I_{e})}
\le C_{\infty}\,\|u\|_{H^{1}(\mathcal G)}
\qquad\forall\,u\in H^{1}(\mathcal G).
\]

We close this subsection with a few basic examples of periodic quantum graphs
that will serve as model geometries.

\begin{Exam}\label{ex-2.1}
	Consider the graph whose vertices are indexed by $\mathbb Z$ and whose edges
	connect consecutive integers. Each edge is identified with an interval
	$[0,\ell]$ of fixed length $\ell>0$. The action of $\mathbb Z$ is generated by
	the shift $T^{1}$ sending the vertex $n$ to $n+1$ and translating each edge by
	one unit. A fundamental cell $\mathcal K$ can be chosen as a single edge
	together with its two endpoints. A schematic picture is given in
	Figure~\ref{fig1}.
	
	\begin{figure}[th]
		\centering
		\begin{tikzpicture}[scale=1]
			\foreach \x in {-2,...,2}
			\fill (\x,0) circle (2pt);
			\foreach \x in {-2,...,1}
			\draw (\x,0)--(\x+1,0);
			\draw[dashed] (-2.5,0)--(-2,0);
			\draw[dashed] (2,0)--(2.5,0);
			\draw[rounded corners] (-0.5,-0.35) rectangle (0.5,0.35);
			\node at (0,-0.8) {$\mathcal K$};
		\end{tikzpicture}
		\caption{A one-dimensional periodic quantum graph (chain) with fundamental
			cell $\mathcal K$.}
		\label{fig1}
	\end{figure}
\end{Exam}

\begin{Exam}\label{ex-2.2}
	Let $\mathcal G$ be the standard square lattice in the plane: vertices are
	points with integer coordinates $(m,n)\in\mathbb Z^{2}$, and each vertex is
	connected to its four nearest neighbours by edges of unit length. The action of
	$\mathbb Z^{2}$ is given by horizontal and vertical translations. A convenient
	fundamental cell $\mathcal K$ is the square with vertices
	$(0,0),(1,0),(1,1),(0,1)$. A schematic picture of a finite portion of the
	lattice is shown in Figure~\ref{fig2}.
	
	\begin{figure}[th]
		\centering
		\begin{tikzpicture}[scale=0.8]
			\foreach \x in {0,...,3}{
				\foreach \y in {0,...,3}{
					\fill (\x,\y) circle (2pt);
				}
			}
			\foreach \x in {0,...,2}{
				\foreach \y in {0,...,3}{
					\draw (\x,\y)--(\x+1,\y);
				}
			}
			\foreach \x in {0,...,3}{
				\foreach \y in {0,...,2}{
					\draw (\x,\y)--(\x,\y+1);
				}
			}
			\draw[rounded corners,thick] (0,-0.2) rectangle (1.2,1.2);
			\node at (0.6,-0.6) {$\mathcal K$};
		\end{tikzpicture}
		\caption{A finite portion of the two-dimensional square lattice and a
			fundamental cell $\mathcal K$.}
		\label{fig2}
	\end{figure}
\end{Exam}

\begin{Exam}\label{ex-2.3}
	A further example of a periodic quantum graph with nontrivial topology is the
	ladder graph, formed by two parallel chains connected by rungs. Vertices are
	$\{(n,0),(n,1): n\in\mathbb Z\}$, with horizontal edges connecting
	$(n,j)$ to $(n+1,j)$ for $j=0,1$ and vertical edges connecting $(n,0)$ to
	$(n,1)$ for all $n\in\mathbb Z$. The action of $\mathbb Z$ is generated by
	translation in the $n$-direction, and a fundamental cell $\mathcal K$ is given
	by two horizontal edges and one vertical rung, see Figure~\ref{fig3}.
	
	\begin{figure}[th]
		\centering
		\begin{tikzpicture}[scale=1]
			\foreach \x in {-1,...,1}{
				\fill (\x,0) circle (2pt);
				\fill (\x,1) circle (2pt);
			}
			\foreach \x in {-1,0}{
				\draw (\x,0)--(\x+1,0);
				\draw (\x,1)--(\x+1,1);
			}
			\foreach \x in {-1,...,1}{
				\draw (\x,0)--(\x,1);
			}
			\draw[dashed] (-1.5,0)--(-1,0);
			\draw[dashed] (-1.5,1)--(-1,1);
			\draw[dashed] (1,0)--(1.5,0);
			\draw[dashed] (1,1)--(1.5,1);
			\draw[rounded corners,thick] (-0.5,-0.3) rectangle (0.5,1.3);
			\node at (0,-0.7) {$\mathcal K$};
		\end{tikzpicture}
		\caption{A periodic ladder graph with fundamental cell $\mathcal K$.}
		\label{fig3}
	\end{figure}
\end{Exam}

\subsection{The Dirac operator with Kirchhoff-type conditions}

The expression of the Dirac operator on a metric graph given in \eqref{eq-1.1}
is purely formal, in the sense that it specifies the action of the operator only
in the interior of each edge, where the derivative $\frac{d}{dx}$ is well
defined. In order to obtain a self-adjoint operator on $L^{2}(\mathcal G,\mathbb C^{2})$
one has to prescribe suitable vertex conditions.

As in the case of the Laplacian in Schr\"odinger equations, there is a large
class of self-adjoint realizations of the Dirac operator on a quantum graph,
described in terms of boundary conditions at the vertices; see, for example,
\cite{MR3013208,MR3934110,MR4219183,MR1050469}. In the present work we focus on
the Kirchhoff-type vertex conditions, which correspond to the ``free'' case for
the Dirac operator.

\begin{Def}\label{def-2.1}
Let $\mathcal G$ be a quantum graph and let $m,c>0$. The Dirac operator with
Kirchhoff-type vertex conditions is the operator
\[
\mathcal D:L^{2}(\mathcal G,\mathbb C^{2})\to L^{2}(\mathcal G,\mathbb C^{2})
\]
with action on each edge $I_{e}$ given by
\begin{equation}\label{eq-2.1}
	\mathcal D|_{I_{e}}u
	=\mathcal D_{e}u_{e}
	=-ic\,\sigma_{1}u_{e}' + mc^{2}\sigma_{3}u_{e},
	\qquad e\in\mathcal E,
\end{equation}
where $\sigma_{1},\sigma_{3}$ are the Pauli matrices, and with domain
\begin{equation}\label{eq-2.2}
	\operatorname{dom}(\mathcal D)
	=\Bigl\{u\in H^{1}(\mathcal G,\mathbb C^{2}):
	u\ \text{satisfies \eqref{eq-2.3} and \eqref{eq-2.4} at every vertex}\Bigr\},
\end{equation}
where, for every vertex $v\in\mathcal V$,
\begin{equation}\label{eq-2.3}
	u_{e}^{1}(v)=u_{f}^{1}(v)
	\quad\text{for all edges }e,f\text{ incident at }v,
\end{equation}
and
\begin{equation}\label{eq-2.4}
	\sum_{e\succ v} u_{e}^{2}(v)_{\pm}=0.
\end{equation}
Here $e\succ v$ indicates that the edge $e$ is incident at the vertex $v$, and
$u_{e}^{2}(v)_{\pm}$ stands for $u_{e}^{2}(0)$ or $-u_{e}^{2}(\ell_{e})$
according to whether the coordinate $x_{e}$ takes the value $0$ or $\ell_{e}$
at $v$.
\end{Def}

\begin{remark}\label{rem-2.2}
The operator $\mathcal D$ depends on the parameters $m$ and $c$, which
represent the mass of the particle and the speed of light, respectively. In the
sequel we suppress this dependence from the notation and simply write
$\mathcal D$. Moreover, in the periodic setting,
the translations $T^{k}$ act unitarily on $L^{2}(\mathcal G,\mathbb C^{2})$ and
\[
\mathcal D(T^{k}u)=T^{k}(\mathcal Du)
\qquad\forall\,k\in\mathbb Z^{d},\ u\in\operatorname{dom}(\mathcal D),
\]
so $\mathcal D$ is $\mathbb Z^{d}$–equivariant.
\end{remark}

We collect here the basic spectral properties of $\mathcal D$ that will be
used in the variational analysis.

\begin{Pro}\label{pro2.1}
Let $\mathcal G$ be a connected metric graph carrying a free, cocompact action
of $\mathbb Z^{d}$ by graph automorphisms as above, and let $\mathcal D$ be the
Dirac operator with Kirchhoff-type vertex conditions introduced in
Definition~\ref{def-2.1}. Then $\mathcal D$ is self-adjoint on
$L^{2}(\mathcal G,\mathbb C^{2})$. Moreover
\begin{equation}\label{eq-2.5}
	\sigma(\mathcal D)\subset(-\infty,-mc^{2}]\cup[mc^{2},+\infty),
\end{equation}
so that $0$ lies in a spectral gap of width $2mc^{2}$.
\end{Pro}

\begin{proof}
Self-adjointness of $\mathcal D$ with Kirchhoff-type vertex conditions is a
special case of the general theory of self-adjoint realizations of the Dirac
operator on metric graphs. For completeness we briefly recall the
main points.

On each edge $I_{e}$ the differential expression $\mathcal D_{e}$ is formally symmetric.
The maximal Dirac operator $\mathcal D_{\max}$ on $\mathcal G$ acts as in
\eqref{eq-2.1} on the domain $H^{1}(\mathcal G,\mathbb C^{2})$, and an
integration by parts shows that for $u,v\in\operatorname{dom}(\mathcal D_{\max})$
one has
\[
\langle \mathcal D_{\max}u,v\rangle_{L^{2}}
-\langle u,\mathcal D_{\max}v\rangle_{L^{2}}
= ic\sum_{v\in\mathcal V}\sum_{e\succ v}
\Bigl(u_{e}^{2}(v)_{\pm}\,\overline{v_{e}^{1}(v)}
-u_{e}^{1}(v)\,\overline{v_{e}^{2}(v)_{\pm}}\Bigr).
\]
The vertex conditions \eqref{eq-2.3}--\eqref{eq-2.4} select at each vertex a maximal
isotropic subspace for this boundary form, hence they define a self-adjoint
restriction of $\mathcal D_{\max}$ (see, e.g., \cite{MR1050469,MR3013208}).

We turn to the spectral gap. On each edge $I_{e}$ one computes the square of the
differential expression,
\[
\mathcal D_{e}^{2}
=\bigl(-ic\,\sigma_{1}\partial_{x_{e}}+mc^{2}\sigma_{3}\bigr)^{2}
=-c^{2}\partial_{x_{e}}^{2}+m^{2}c^{4}I_{2},
\]
since $\sigma_{1}^{2}=\sigma_{3}^{2}=I_{2}$ and $\sigma_{1}\sigma_{3}+\sigma_{3}\sigma_{1}=0$.
Moreover, using $\operatorname{dom}(\mathcal D^{2})
=\{u\in\operatorname{dom}(\mathcal D):\ \mathcal Du\in\operatorname{dom}(\mathcal D)\}$,
one checks that the two components decouple at the level of $\mathcal D^{2}$.
More precisely, for every vertex $v\in\mathcal V$ one has
\[
u^{1}\ \text{is continuous at }v,
\qquad
\sum_{e\succ v} (u^{1}_{e})'(v)_{\pm}=0,
\]
and
\[
\sum_{e\succ v} u^{2}_{e}(v)_{\pm}=0,
\qquad
(u^{2}_{e})'(v)_{\pm}=(u^{2}_{f})'(v)_{\pm}\ \text{for all }e,f\succ v.
\]
Thus $\mathcal D^{2}$ can be written as the diagonal operator
\[
\mathcal D^{2}
=
\begin{pmatrix}
c^{2}(-\Delta_{\mathrm K})+m^{2}c^{4} & 0\\
0 & c^{2}(-\Delta_{\mathrm{AK}})+m^{2}c^{4}
\end{pmatrix},
\]
where $-\Delta_{\mathrm K}$ is the Kirchhoff Laplacian, and $-\Delta_{\mathrm{AK}}$ is the anti-Kirchhoff Laplacian.
By integration by parts the boundary terms vanish under either set of vertex conditions, and one obtains
$\langle -\Delta_{\mathrm K}f,f\rangle_{L^{2}}=\|f'\|_{L^{2}}^{2}\ge0$ and
$\langle -\Delta_{\mathrm{AK}}f,f\rangle_{L^{2}}=\|f'\|_{L^{2}}^{2}\ge0$ on their
respective domains. Consequently,
\[
\sigma(\mathcal D^{2}) \subset [m^{2}c^{4},+\infty).
\]
By the spectral mapping theorem for self-adjoint operators,
\[
\sigma(\mathcal D)
= \{\lambda\in\mathbb R:\ \lambda^{2}\in\sigma(\mathcal D^{2})\}
\subset (-\infty,-mc^{2}]\cup[mc^{2},+\infty),
\]
which is \eqref{eq-2.5}. In particular $(-mc^{2},mc^{2})$ contains no spectral
points of $\mathcal D$, so that $0$ lies in a spectral gap of width $2mc^{2}$.
\end{proof}

\begin{remark}\label{rem-2.3}
In the $\mathbb Z^{d}$-periodic setting described above, the operator $\mathcal D$
is $\mathbb Z^{d}$--equivariant and therefore admits a Floquet--Bloch decomposition.
As a consequence, $\sigma(\mathcal D)$ has a band-gap structure contained in the two
half-lines $(-\infty,-mc^{2}]$ and $[mc^{2},+\infty)$; see, for instance,
\cite{KuchmentQG2}. In particular, for suitable periodic geometries of $\mathcal G$
there may exist internal spectral gaps in $[mc^{2},+\infty)$ and in $(-\infty,-mc^{2}]$.

On the other hand, there are periodic quantum graphs for which the Kirchhoff Laplacian
has no gaps at all and $\sigma(-\Delta_{\mathrm K})=[0,+\infty)$; see, for example,
the rectangular graph superlattices discussed in \cite{ExnerContact96,ExnerGawlista96}.
In several such geometries the inclusion in \eqref{eq-2.5} is sharp, namely one has
\[
\sigma(\mathcal D)=(-\infty,-mc^{2}]\cup[mc^{2},+\infty).
\]
\end{remark}

\subsection{The associated quadratic form}

A convenient way to describe the form domain of the Dirac operator
$\mathcal D$ is provided by interpolation theory
(see, e.g., \cite{MR2424078,MR2715493}). We only recall the basic facts
needed in this paper and refer to \cite{MR3934110} for further details.

We set
\begin{equation}\label{eq-2.6}
  Y = \bigl[L^{2}(\mathcal G,\mathbb C^{2}),\ \operatorname{dom}(\mathcal D)\bigr]_{1/2},
\end{equation}
that is, $Y$ is the interpolation space of order $1/2$ between
$L^{2}(\mathcal G,\mathbb C^{2})$ and the domain of the Dirac operator.
Since $\operatorname{dom}(\mathcal D)$ is a closed subspace of
$H^{1}(\mathcal G,\mathbb C^{2})$ and
\[
H^{1/2}(\mathcal G,\mathbb C^{2})
= \bigl[L^{2}(\mathcal G,\mathbb C^{2}),\ H^{1}(\mathcal G,\mathbb C^{2})\bigr]_{1/2},
\]
it follows from the general properties of interpolation of closed subspaces that
$Y$ can be identified with a closed subspace of $H^{1/2}(\mathcal G,\mathbb C^{2})$,
endowed with the norm induced by $H^{1/2}(\mathcal G,\mathbb C^{2})$.

By Sobolev embeddings on one-dimensional intervals and the periodic structure of
$\mathcal G$ (in particular, the uniform bound on edge lengths coming from the
finiteness of the fundamental cell $\mathcal K$), we obtain
\begin{equation}\label{eq-2.7}
  Y \hookrightarrow L^{p}(\mathcal G,\mathbb C^{2})
  \qquad\text{for all }p\in[2,+\infty),
\end{equation}
and, moreover, the embedding
\[
Y \hookrightarrow L^{p}(\mathcal K,\mathbb C^{2})
\]
is compact for every $p\in[2,+\infty)$, thanks to the compactness of
$\mathcal K$ and the Rellich theorem on finite unions of intervals.

On the other hand, one has (see \cite{MR3934110})
\begin{equation}\label{eq-2.8}
  \operatorname{dom}\bigl(\mathcal Q_{\mathcal D}\bigr) = Y,
\end{equation}
where $\mathcal Q_{\mathcal D}$ denotes the (closed) quadratic form associated with
$\mathcal D$ via spectral calculus. This identification of the form domain with
the interpolation space is a key point in the variational analysis below.

For later use we fix the notation
\[
\mathcal Q_{\mathcal D}(u)
= \frac{1}{2}\int_{\mathcal G}\langle u,\mathcal D u\rangle\,dx,
\qquad
\mathcal Q_{\mathcal D}(u,v)
= \frac{1}{2}\int_{\mathcal G}\langle u,\mathcal D v\rangle\,dx,
\]
for $u,v\in Y$, where $\langle\cdot,\cdot\rangle$ denotes the standard
sesquilinear scalar product on $\mathbb C^{2}$. Whenever $u$ and $v$ are
sufficiently regular, in particular when $u,v\in\operatorname{dom}(\mathcal D)$,
these expressions coincide with the abstract quadratic form and the associated
bilinear form defined by $\mathcal Q_{\mathcal D}$.

\section{Variational setting and abstract critical point theory}

\subsection{Variational setting}
We now introduce the variational functional associated with \eqref{eq-1.4}.
For $u\in Y$ we consider
\begin{equation}\label{eq-3.1}
\Phi(u)
= \frac{1}{2}\mathcal Q_{\mathcal D}(u,u)
+ \frac{\omega}{2}\int_{\mathcal G}|u|^{2}\,dx
- \int_{\mathcal G} F(x,u)\,dx .
\end{equation}
If $u\in\operatorname{dom}(\mathcal D)$, then $\mathcal Q_{\mathcal D}(u,u)=(\mathcal Du,u)_{L^{2}}$,
and \eqref{eq-3.1} agrees with the formal expression
$\frac12\int_{\mathcal G}\langle u,(\mathcal D+\omega)u\rangle\,dx-\int_{\mathcal G}F(x,u)\,dx$.
Under assumptions $(F_{0})$--$(F_{5})$ one checks that $\Phi\in C^{2}(Y,\mathbb R)$.

Recall from Proposition~\ref{pro2.1} that
\begin{equation}\label{eq-3.2}
	\sigma(\mathcal D)\subset(-\infty,-mc^{2}]\cup[mc^{2},+\infty),
\end{equation}
so that $0$ lies in a spectral gap. By the spectral theorem, the Hilbert space
$Y$ decomposes as the orthogonal sum of the positive and negative spectral
subspaces of $\mathcal D$,
\[
Y = Y^{+}\oplus Y^{-},
\]
where $Y^{\pm}$ are the ranges of the spectral projectors $P^{\pm}$ associated
with $(0,+\infty)$ and $(-\infty,0)$, respectively. Thus every $u\in Y$ can be
written uniquely as
\[
u = u^{+}+u^{-},\qquad
u^{\pm} = P^{\pm}u .
\]
It is convenient to equip $Y$ with the equivalent norm
\[
\|u\|^{2}
= \bigl\|\lvert\mathcal D\rvert^{1/2}u\bigr\|_{L^{2}(\mathcal G,\mathbb C^{2})}^{2}
= \bigl(\lvert\mathcal D\rvert u,u\bigr)_{L^{2}(\mathcal G,\mathbb C^{2})},
\qquad u\in Y.
\]

The next estimate will be used frequently.

\begin{lemma}\label{lem-3.1}
	For every $u\in Y$ one has
	\begin{equation}\label{eq-3.3}
		mc^{2}\,\|u\|_{L^{2}(\mathcal G,\mathbb C^{2})}^{2}
		\le \|u\|^{2}.
	\end{equation}
\end{lemma}

\begin{proof}
	By \eqref{eq-3.2} the spectrum of $\lvert\mathcal D\rvert$ is contained in
	$[mc^{2},+\infty)$, hence $\lvert\mathcal D\rvert\ge mc^{2}I$ as a self-adjoint
	operator on $L^{2}(\mathcal G,\mathbb C^{2})$. Therefore, for every $u\in Y$,
	\[
	\|u\|^{2}
	= \bigl(\lvert\mathcal D\rvert u,u\bigr)_{L^{2}(\mathcal G,\mathbb C^{2})}
	\ge mc^{2}\,\|u\|_{L^{2}(\mathcal G,\mathbb C^{2})}^{2},
	\]
	which is exactly \eqref{eq-3.3}.
\end{proof}

We now relate bound states of \eqref{eq-1.4} to critical points of $\Phi$.

\begin{Pro}\label{pro-3.1}
	A spinor $u$ is a bound state of frequency $\omega$ of the NLDE \eqref{eq-1.4}
	if and only if $u$ is a critical point of $\Phi$.
\end{Pro}

\begin{proof}
Assume first that $u\in\operatorname{dom}(\mathcal D)$ is a bound state of frequency $\omega$,
namely
\[
(\mathcal D+\omega)u = F_{u}(x,u)\quad\text{in }L^{2}(\mathcal G,\mathbb C^{2}).
\]
Then $u\in Y$ and for every $\varphi\in Y$ we have
\[
\Phi'(u)[\varphi]
=\mathcal Q_{\mathcal D}(u,\varphi)
+\omega\int_{\mathcal G}\langle u,\varphi\rangle\,dx
-\int_{\mathcal G}\langle F_{u}(x,u),\varphi\rangle\,dx.
\]
Since $u\in\operatorname{dom}(\mathcal D)$, one has $\mathcal Q_{\mathcal D}(u,\varphi)=(\mathcal Du,\varphi)_{L^{2}}$,
hence
\[
\Phi'(u)[\varphi]=\bigl((\mathcal D+\omega)u-F_{u}(x,u),\varphi\bigr)_{L^{2}}=0
\qquad\forall\,\varphi\in Y,
\]
so $u$ is a critical point of $\Phi$.

Conversely, let $u\in Y$ be a critical point of $\Phi$. Then
\begin{equation}\label{eq-3.4}
\mathcal Q_{\mathcal D}(u,\varphi)
+\omega\int_{\mathcal G}\langle u,\varphi\rangle\,dx
=\int_{\mathcal G}\langle F_{u}(x,u),\varphi\rangle\,dx
\qquad\forall\,\varphi\in Y.
\end{equation}
In particular, \eqref{eq-3.4} holds for all $\varphi\in\operatorname{dom}(\mathcal D)\subset Y$.
For such $\varphi$ one has $\mathcal Q_{\mathcal D}(u,\varphi)=(u,\mathcal D\varphi)_{L^{2}}$, hence
\[
\bigl(u,(\mathcal D+\omega)\varphi\bigr)_{L^{2}}
=\bigl(F_{u}(x,u),\varphi\bigr)_{L^{2}}
\qquad\forall\,\varphi\in\operatorname{dom}(\mathcal D).
\]
Therefore $u\in\operatorname{dom}((\mathcal D+\omega)^{*})=\operatorname{dom}(\mathcal D)$ and
\[
(\mathcal D+\omega)u=F_{u}(x,u)\quad\text{in }L^{2}(\mathcal G,\mathbb C^{2}).
\]
Thus $u$ is a bound state of frequency $\omega$ of \eqref{eq-1.4}.
\end{proof}

In summary, using the spectral decomposition $Y=Y^{+}\oplus Y^{-}$ and the norm
$\|u\|^{2}=\bigl(\lvert\mathcal D\rvert u,u\bigr)_{L^{2}(\mathcal G,\mathbb C^{2})}$,
we can rewrite the action functional \eqref{eq-3.1} in the form
\begin{equation}\label{eq-3.6}
	\Phi(u)
	= \frac{1}{2}\bigl(\|u^{+}\|^{2}-\|u^{-}\|^{2}\bigr)
	+ \frac{\omega}{2}\int_{\mathcal G}|u|^{2}\,dx
	- \Psi(u),
\end{equation}
where
\[
\Psi(u) = \int_{\mathcal G} F(x,u)\,dx.
\]
Here $u^{\pm}=P^{\pm}u$ are the components of $u$ in the positive and negative
spectral subspaces $Y^{\pm}$ of $\mathcal D$, and
$\|u^{\pm}\|^{2}=\bigl(\lvert\mathcal D\rvert u^{\pm},u^{\pm}\bigr)_{L^{2}(\mathcal G,\mathbb C^{2})}$.

\subsection{Critical point theorems}

We recall the abstract critical point framework that we shall use, following
\cite{MR2255874}. Let $Z$ be a Banach space with a topological direct sum
decomposition
\[
Z = M \oplus N
\]
and corresponding continuous projections $P_{M},P_{N}$ onto $M$ and $N$,
respectively. For $\Phi\in C^{1}(Z,\mathbb R)$ and $a,b\in\mathbb R$ we set
\[
\Phi_{a} = \{u\in Z:\ \Phi(u)\ge a\},\qquad
\Phi^{b} = \{u\in Z:\ \Phi(u)\le b\},\qquad
\Phi_{a}^{b} = \Phi_{a}\cap\Phi^{b}.
\]

\begin{Def}\label{def-3.1}
	A sequence $(u_{n})\subset Z$ is called a $(C)_{c}$–sequence if
	\[
	\Phi(u_{n})\to c
	\quad\text{and}\quad
	\bigl(1+\|u_{n}\|\bigr)\,\|\Phi'(u_{n})\|_{Z^{*}}\to 0.
	\]
	We say that $\Phi$ satisfies the $(C)_{c}$–condition if every $(C)_{c}$–sequence
	has a convergent subsequence.
\end{Def}

\begin{Def}\label{def-3.2}
	A set $\mathcal A\subset Z$ is called a $(C)_{c}$–attractor if for every
	$\varepsilon,\delta>0$ and every $(C)_{c}$–sequence $(u_{n})$ there exists
	$n_{0}$ such that
	\[
	u_{n}\in U_{\varepsilon}\bigl(\mathcal A\cap\Phi_{c-\delta}^{c+\delta}\bigr)
	\qquad\text{for all }n\ge n_{0},
	\]
	where $U_{\varepsilon}(B)$ denotes the $\varepsilon$–neighbourhood of a set
	$B\subset Z$. Given an interval $I\subset\mathbb R$, we say that $\mathcal A$
	is a $(C)_{I}$–attractor if it is a $(C)_{c}$–attractor for every $c\in I$.
\end{Def}

From now on we assume that $M$ is separable and reflexive, and we fix a countable dense
subset $\mathcal S\subset M^{*}$. For each $s\in\mathcal S$ we define a seminorm
on $Z$ by
\[
p_{s}(u) = |s(x)| + \|y\|
\quad\text{for }u=x+y\in M\oplus N.
\]
Here $s\in M^{*}$ acts only on the $M$--component $x=P_{M}u$.
We denote by $\mathcal T_{\mathcal S}$ the locally convex topology on $Z$
generated by the family $\{p_{s}: s\in\mathcal S\}$, and by $w^{*}$ the weak*
topology on $Z^{*}$.

We shall use the following structural assumptions on $\Phi$:

\begin{itemize}
	\item[$(\Phi_{0})$]
	For every $c\in\mathbb R$, the set $\Phi_{c}$ is
	$\mathcal T_{\mathcal S}$–closed, and
	\[
	\Phi' : (\Phi_{c},\mathcal T_{\mathcal S}) \to (Z^{*},w^{*})
	\]
	is continuous.
	
	\item[$(\Phi_{1})$]
	For every $c>0$ there exists $\zeta>0$ such that
	\[
	\|u\| < \zeta\,\|P_{N}u\|
	\qquad\text{for all }u\in\Phi_{c}.
	\]
	
	\item[$(\Phi_{2})$]
	There exists $\rho>0$ such that
	\[
	\kappa= \inf_{u\in S_{\rho}^{N}} \Phi(u) > 0,
	\]
	where
	\[
	S_{\rho}^{N} = \{u\in N:\ \|u\|=\rho\}.
	\]
	
	\item[$(\Phi_{3})$]
	There exist a finite-dimensional subspace $N_{0}\subset N$ and a number
	$R>\rho$ such that, denoting
	\[
	E_{0} = M\oplus N_{0}, \qquad
	B_{0} = \{u\in E_{0}:\ \|u\|\le R\},
	\]
	we have $b' := \sup\Phi(E_{0})<\infty$ and
	\[
	\sup\Phi(E_{0}\setminus B_{0})
	< \inf\Phi\bigl(B_{\rho}^{N}\bigr),
	\]
	where $B_{\rho}^{N}=\{u\in N:\ \|u\|\le\rho\}$.
	
	\item[$(\Phi_{4})$]
	There exist an increasing sequence of finite-dimensional subspaces
	$N_{n}\subset N$ and a sequence $(R_{n})$ of positive numbers such that,
	setting
	\[
	E_{n} = M\oplus N_{n},\qquad
	  B_{n} = \{u\in E_{n}:\ \|u\|\le R_{n}\},
	\]
	one has $\sup\Phi(E_{n})<\infty$ and
	\[
	\sup\Phi(E_{n}\setminus B_{n})
	< \inf\Phi\bigl(B_{\rho}^{N}\bigr)
	\qquad\text{for all }n.
	\]
	
	\item[$(\Phi_{5})$]
	One of the following alternatives holds:
	\begin{itemize}
		\item[(i)] For every interval $I\subset(0,\infty)$ there exists a
		$(C)_{I}$–attractor $\mathcal A$ such that $P_{N}\mathcal A$ is bounded and
		\[
		\inf\Bigl\{\|P_{N}(u-v)\|:\ u,v\in\mathcal A,\ P_{N}(u-v)\ne 0\Bigr\} > 0;
		\]
		\item[(ii)] $\Phi$ satisfies the $(C)_{c}$–condition for all $c>0$.
	\end{itemize}
\end{itemize}

The following result is a generalized linking theorem; it is a special case of
\cite{MR2255874}.

\begin{theorem}\label{theo-3.1}
	Assume that $(\Phi_{0})$–$(\Phi_{2})$ hold. Suppose there exist numbers
	$R>\rho>0$ and an element $e\in N$ with $\|e\|=1$ such that
	\[
	\sup\Phi(\partial Q)\le\kappa,
	\]
	where
	\[
	Q = \{u = x + te:\ x\in M,\ t\ge 0,\ \|u\|<R\}.
	\]
	Then $\Phi$ has a $(C)_{c}$–sequence with
	\[
	\kappa \le c \le \sup\Phi(Q).
	\]
\end{theorem}

\begin{theorem}\label{theo-3.2}
	Assume that $\Phi$ is even and $\Phi(0)=0$, and that
	$(\Phi_{0})$–$(\Phi_{5})$ are satisfied. Then $\Phi$ possesses an unbounded
	sequence of positive critical values.
\end{theorem}

In addition, we recall a convenient criterion for verifying $(\Phi_{0})$.

\begin{theorem}\label{theo-3.3}
	Let $Z = M\oplus N$ as above and suppose that $\Phi\in C^{1}(Z,\mathbb R)$ is of
	the form
	\[
	\Phi(u) = \frac{1}{2}\bigl(\|y\|^{2}-\|x\|^{2}\bigr) - \Psi(u)
	\quad\text{for }u=x+y\in M\oplus N,
	\]
	where $\Psi\in C^{1}(Z,\mathbb R)$ satisfies:
	\begin{itemize}
		\item[(i)] $\Psi$ is bounded from below;
		\item[(ii)] $\Psi:(Z,\mathcal T_{\mathcal S})\to\mathbb R$ is sequentially
		lower semicontinuous, that is,
		\[
		u_{n}\to u\ \text{in }(Z,\mathcal T_{\mathcal S})
		\ \to\ 
		\Psi(u)\le\liminf_{n\to\infty}\Psi(u_{n});
		\]
		\item[(iii)] $\Psi':(Z,\mathcal T_{\mathcal S})\to(Z^{*},w^{*})$ is
		sequentially continuous;
		\item[(iv)] the map $v:Z\to\mathbb R$, $v(u)=\|u\|^{2}$, is of class $C^{1}$
		and $v':(Z,\mathcal T_{\mathcal S})\to(Z^{*},w^{*})$ is sequentially
		continuous.
	\end{itemize}
	Then $\Phi$ satisfies $(\Phi_{0})$.
\end{theorem}

In order to handle the nonlinear terms and the lack of compactness caused by the
periodic structure of the graph, we establish the following Brezis--Lieb type lemma and a concentration--compactness principle on periodic
quantum graphs.

\begin{lemma}\label{lem-3.2}
	Let $\mathcal G$ be a noncompact metric graph, and let
	$1<p<\infty$. Suppose $\{u_{n}\}\subset L^{p}(\mathcal G)$ and $u\in L^{p}(\mathcal G)$
	satisfy
	\begin{enumerate}
		\item[\textnormal{(i)}] $u_{n}(x)\to u(x)$ almost everywhere on $\mathcal G$,
		\item[\textnormal{(ii)}] $\sup_{n}\|u_{n}\|_{L^{p}(\mathcal G)}<\infty$.
	\end{enumerate}
	Then, as $n\to\infty$,
	\[
	\|u_{n}\|_{L^{p}(\mathcal G)}^{p}
	= \|u_{n}-u\|_{L^{p}(\mathcal G)}^{p}
	+ \|u\|_{L^{p}(\mathcal G)}^{p}
	+ o(1).
	\]
\end{lemma}

\begin{proof}
	The metric graph $(\mathcal G,dx)$ is a $\sigma$–finite measure space, since
	it can be written as a countable union of edges $I_{e}$ with the standard
	one–dimensional Lebesgue measure. The classical Brezis--Lieb lemma holds on
	arbitrary $\sigma$–finite measure spaces (see, for instance,
	\cite{MR1400007}). Applying that result with $u_{n}$ and $u$ on
	$(\mathcal G,dx)$ yields exactly the desired decomposition.
\end{proof}

\begin{lemma}\label{Lem-3.3}
	Let $\mathcal G$ be a connected noncompact metric graph admitting a free,
	cocompact action $\{T^{a}\}_{a\in\mathbb Z^{d}}$ of $\mathbb Z^{d}$ by graph
	isometries, and let $\mathcal K\subset\mathcal G$ be a fixed fundamental cell.
	Let $\{\Psi_{n}\}\subset H^{1}(\mathcal G)$ satisfy
	\begin{itemize}
		\item $\sup_{n}\|\Psi_{n}\|_{H^{1}(\mathcal G)}<\infty$,
		\item $\|\Psi_{n}\|_{L^{2}(\mathcal G)}^{2}\to m>0$.
	\end{itemize}
	Then, up to a subsequence, exactly one of the following mutually exclusive
	alternatives holds:
	\begin{enumerate}
		\item[\textnormal{(i)}] \textnormal{Vanishing.} For every $R>0$,
		\[
		\lim_{n\to\infty}\sup_{x\in\mathcal G}
		\int_{B_{R}(x)}|\Psi_{n}(y)|^{2}\,dy = 0.
		\]
		In particular, $\Psi_{n}\to 0$ strongly in $L^{p}(\mathcal G)$ for all
		$p\in(2,\infty)$.
		
		\item[\textnormal{(ii)}] \textnormal{Dichotomy.}
		There exists $\alpha\in(0,m)$ and sequences $\{R_{n}\},\{S_{n}\}\subset H^{1}(\mathcal G)$
		such that:
		\begin{itemize}
			\item $\Psi_{n} = R_{n} + S_{n} + o_{L^{2}}(1)$, i.e.
			$\|\Psi_{n}-R_{n}-S_{n}\|_{L^{2}(\mathcal G)}\to 0$,
			\item $\|R_{n}\|_{L^{2}(\mathcal G)}^{2}\to\alpha$,
			$\ \|S_{n}\|_{L^{2}(\mathcal G)}^{2}\to m-\alpha$,
			\item $\inf_{n}\operatorname{dist}\bigl(\operatorname{supp}R_{n},
			\operatorname{supp}S_{n}\bigr)>0$,
			\item for all $p\in(2,\infty)$,
			\[
			\|\Psi_{n}\|_{L^{p}(\mathcal G)}^{p}
			= \|R_{n}\|_{L^{p}(\mathcal G)}^{p}
			+ \|S_{n}\|_{L^{p}(\mathcal G)}^{p}
			+ o(1),
			\]
			\[
			\|\Psi_{n}'\|_{L^{2}(\mathcal G)}^{2}
			= \|R_{n}'\|_{L^{2}(\mathcal G)}^{2}
			+ \|S_{n}'\|_{L^{2}(\mathcal G)}^{2}
			+ o(1).
			\]
		\end{itemize}
		
		\item[\textnormal{(iii)}] \textnormal{Compactness modulo translations.}
		There exist a sequence $\{a_{n}\}\subset\mathbb Z^{d}$ and a function
		$\Psi\in H^{1}(\mathcal G)$ such that, up to a subsequence,
		\[
		V_{n} = T^{-a_{n}}\Psi_{n} \rightharpoonup \Psi
		\text{ in }H^{1}(\mathcal G),
		\quad
		V_{n}\to\Psi\text{ in }L^{p}_{\mathrm{loc}}(\mathcal G)
		\text{ for all }p\in[2,\infty),
		\]
		and
		\[
		\|\Psi\|_{L^{2}(\mathcal G)}^{2} = m.
		\]
	\end{enumerate}
\end{lemma}

\begin{proof}
	For every $R>0$ set
	\[
	\rho_{n}(R)
	= \sup_{x\in\mathcal G}\int_{B_{R}(x)}|\Psi_{n}(y)|^{2}\,dy.
	\]
	Each $\rho_{n}$ is nondecreasing in $R$ and satisfies
	\[
	0\le\rho_{n}(R)\le\|\Psi_{n}\|_{L^{2}(\mathcal G)}^{2}\to m
	\quad\text{as }n\to\infty.
	\]
	Define
	\[
	\tau
	= \lim_{R\to\infty}\,\liminf_{n\to\infty}\rho_{n}(R)\in[0,m].
	\]
	
	\medskip \noindent\textbf{Step 1: $\tau=0$.}
Since the map $R\mapsto \liminf\limits_{n\to\infty}\rho_n(R)$ is nondecreasing and
\[
\tau=\lim_{R\to\infty}\,\liminf_{n\to\infty}\rho_n(R)=0,
\]
we have $\liminf\limits_{n\to\infty}\rho_n(k)=0$ for every integer $k\in\mathbb N$.
By a diagonal argument we may pass to a subsequence (still denoted $\{\Psi_n\}$) such that
\[
\rho_n(k)\to0 \qquad\text{for every }k\in\mathbb N.
\]
Given $R>0$, choose $k\in\mathbb N$ with $k>R$. Since $\rho_n$ is nondecreasing,
\[
0\le\rho_n(R)\le\rho_n(k)\to0,
\]
which gives the vanishing condition in $(i)$.

	To prove the $L^{p}$–convergence, use the periodic decomposition of $\mathcal G$.
	By cocompactness there exists a finite fundamental cell $\mathcal K$ such that
	\[
	\mathcal G = \bigcup_{a\in\mathbb Z^{d}}T^{a}(\mathcal K),
	\]
	where the union is disjoint up to vertices and each $T^{a}(\mathcal K)$ is a finite
	union of edges. Let $R_{0}>\operatorname{diam}(\mathcal K)$. For each $a\in\mathbb Z^{d}$
	choose $x_{a}\in T^{a}(\mathcal K)$; then $T^{a}(\mathcal K)\subset B_{R_{0}}(x_{a})$.
	By vanishing at radius $R_{0}$,
	\[
	\int_{T^{a}(\mathcal K)}|\Psi_{n}|^{2}\,dx
	\le \rho_{n}(R_{0})\xrightarrow[n\to\infty]{}0
	\quad\text{uniformly in }a.
	\]
	Set
	\[
	a_{n,a} = \int_{T^{a}(\mathcal K)}|\Psi_{n}|^{2}\,dx,\qquad a\in\mathbb Z^{d}.
	\]
	Then $\sup_{a}a_{n,a}\to0$ and
	\[
	\sum_{a\in\mathbb Z^{d}}a_{n,a}
	= \|\Psi_{n}\|_{L^{2}(\mathcal G)}^{2}\to m.
	\]
	
	The $L^{p}$–norm can be written as
	\[
	\|\Psi_{n}\|_{L^{p}(\mathcal G)}^{p}
	= \sum_{a\in\mathbb Z^{d}}\int_{T^{a}(\mathcal K)}|\Psi_{n}|^{p}\,dx.
	\]
    On each translate $T^{a}(\mathcal K)$ the one--dimensional
Gagliardo--Nirenberg inequality yields, for every $p\in(2,\infty)$, a constant $C_{p}>0$
depending only on $\mathcal K$ such that
\[
\|\Psi_{n}\|_{L^{p}(T^{a}(\mathcal K))}^{p}
\le C_{p}\,\|\Psi_{n}\|_{L^{2}(T^{a}(\mathcal K))}^{p-2}\,
\|\Psi_{n}\|_{H^{1}(T^{a}(\mathcal K))}^{2}.
\]
Set
\[
a_{n,a}=\|\Psi_{n}\|_{L^{2}(T^{a}(\mathcal K))}^{2},
\qquad
b_{n,a}=\|\Psi_{n}'\|_{L^{2}(T^{a}(\mathcal K))}^{2}.
\]
Then $\|\Psi_{n}\|_{H^{1}(T^{a}(\mathcal K))}^{2}=a_{n,a}+b_{n,a}$ and hence
\[
\int_{T^{a}(\mathcal K)}|\Psi_{n}|^{p}\,dx
\le C_{p}\,a_{n,a}^{\frac{p-2}{2}}(a_{n,a}+b_{n,a})
= C_{p}\Bigl(a_{n,a}^{\frac p2}+a_{n,a}^{\frac p2-1}b_{n,a}\Bigr).
\]
Summing over $a\in\mathbb Z^{d}$ gives
\[
\|\Psi_{n}\|_{L^{p}(\mathcal G)}^{p}
\le C_{p}\sum_{a}a_{n,a}^{\frac p2}
+ C_{p}\sum_{a}a_{n,a}^{\frac p2-1}b_{n,a}.
\]
Using $\sum_{a}a_{n,a}=\|\Psi_{n}\|_{L^{2}(\mathcal G)}^{2}$ and
$\sum_{a}b_{n,a}=\|\Psi_{n}'\|_{L^{2}(\mathcal G)}^{2}$, together with
\[
\sum_{a}a_{n,a}^{\gamma}\le\bigl(\sup_{a}a_{n,a}\bigr)^{\gamma-1}\sum_{a}a_{n,a}
\qquad(\gamma>1),
\]
we obtain (with $\gamma=p/2>1$)
\[
\sum_{a}a_{n,a}^{\frac p2}\le\bigl(\sup_{a}a_{n,a}\bigr)^{\frac p2-1}\|\Psi_{n}\|_{2}^{2},
\qquad
\sum_{a}a_{n,a}^{\frac p2-1}b_{n,a}\le\bigl(\sup_{a}a_{n,a}\bigr)^{\frac p2-1}\|\Psi_{n}'\|_{2}^{2}.
\]
Therefore
\[
\|\Psi_{n}\|_{L^{p}(\mathcal G)}^{p}
\le C_{p}\bigl(\sup_{a}a_{n,a}\bigr)^{\frac p2-1}\bigl(\|\Psi_{n}\|_{2}^{2}+\|\Psi_{n}'\|_{2}^{2}\bigr).
\]
Since $a_{n,a}\le\rho_{n}(R_{0})\to0$ uniformly in $a$ and $\sup_n\|\Psi_n\|_{H^{1}(\mathcal G)}<\infty$,
we conclude that $\Psi_{n}\to0$ in $L^{p}(\mathcal G)$ for all $p\in(2,\infty)$ and
	alternative (i) holds when $\tau=0$.
	
\medskip \noindent\textbf{Step 2: $0<\tau<m$.}
In this case we obtain dichotomy. Fix a decreasing sequence
$\varepsilon_k\rightarrow0$. For each $k\in\mathbb N$ we choose radii
$0<R_{1,k}<R_{2,k}<R_{3,k}$ and an index $n_k$ as follows.

Since
\[
\tau=\lim_{R\to\infty}\,\liminf_{n\to\infty}\rho_n(R)\in(0,m),
\]
we can first choose $R_{1,k}$ so large that
\[
\liminf_{n\to\infty}\rho_n(R_{1,k})>\tau-\varepsilon_k.
\]
Next, using the monotonicity of $\rho_n$ in $R$ and the convergence of
$\liminf_{n}\rho_n(R)$ to $\tau$, we can choose $R_{3,k}>R_{1,k}$ large enough
so that
\[
\liminf_{n\to\infty}\rho_n(R_{3,k})<\tau+\varepsilon_k
\quad\text{and}\quad
\liminf_{n\to\infty}\rho_n(R_{3,k})-\liminf_{n\to\infty}\rho_n(R_{1,k})
<\varepsilon_k.
\]
Passing to a subsequence in $n$ if necessary, we may assume that
\[
\rho_n(R_{1,k})\to\alpha_k,\qquad
\rho_n(R_{3,k})\to\beta_k
\]
with
\[
\tau-\varepsilon_k<\alpha_k\le\beta_k<\tau+\varepsilon_k,\qquad
0<\alpha_k<m.
\]
In particular
\[
0\le\beta_k-\alpha_k<\varepsilon_k.
\]

For each $k$ and each $n$ large enough, we can pick a point $x_{n,k}\in\mathcal G$ such
that
\[
\int_{B_{R_{1,k}}(x_{n,k})}|\Psi_n|^2\,dx
\ge\rho_n(R_{1,k})-\varepsilon_k.
\]
Then
\[
\begin{aligned}
\int_{B_{R_{3,k}}(x_{n,k})\setminus B_{R_{1,k}}(x_{n,k})}|\Psi_n|^2\,dx
&=\int_{B_{R_{3,k}}(x_{n,k})}|\Psi_n|^2\,dx
 -\int_{B_{R_{1,k}}(x_{n,k})}|\Psi_n|^2\,dx\\
&\le\rho_n(R_{3,k})-\bigl(\rho_n(R_{1,k})-\varepsilon_k\bigr).
\end{aligned}
\]
Taking $n\to\infty$ and using the convergence of $\rho_n(R_{j,k})$ we get
\[
\limsup_{n\to\infty}
\int_{B_{R_{3,k}}(x_{n,k})\setminus B_{R_{1,k}}(x_{n,k})}|\Psi_n|^2\,dx
\le\beta_k-\alpha_k+\varepsilon_k<2\varepsilon_k.
\]

Now fix $k$ and choose $R_{2,k}$ such that
\[
R_{1,k}+1<R_{2,k}<R_{3,k}-1.
\]
Choose cut--off functions $\chi_{n,k},\eta_{n,k}\in H^{1}(\mathcal G)\cap L^{\infty}(\mathcal G)$ with
\[
0\le\chi_{n,k}\le1,\quad
\chi_{n,k}\equiv1\ \text{on }B_{R_{1,k}}(x_{n,k}),\quad
\operatorname{supp}\chi_{n,k}\subset B_{R_{2,k}}(x_{n,k}),
\]
\[
0\le\eta_{n,k}\le1,\quad
\eta_{n,k}\equiv0\ \text{on }B_{R_{3,k}-1}(x_{n,k}),\quad
\eta_{n,k}\equiv1\ \text{on }\mathcal G\setminus B_{R_{3,k}}(x_{n,k}),
\]
and
\[
\|\chi_{n,k}'\|_{L^{\infty}}+\|\eta_{n,k}'\|_{L^{\infty}}\le C
\]
for some constant $C>0$ independent of $n,k$.
Define
\[
R_{n,k}=\chi_{n,k}\Psi_n,\quad
S_{n,k}=\eta_{n,k}\Psi_n,\quad
W_{n,k}=\Psi_n-R_{n,k}-S_{n,k}.
\]
Then $\operatorname{dist}(\operatorname{supp}R_{n,k},\operatorname{supp}S_{n,k})\ge (R_{3,k}-1)-R_{2,k}>0$.

Finally, we perform a diagonal extraction: choose an increasing sequence
$n_k$ such that all the above convergences hold along $\{\Psi_{n_k}\}$ and
$\varepsilon_k\to0$. Relabelling
\[
\Psi_{n_k}\mapsto\Psi_n,\quad
R_{n_k,k}\mapsto R_n,\quad
S_{n_k,k}\mapsto S_n,
\]
we obtain sequences $R_n,S_n\in H^1(\mathcal G)$ such that
\[
\Psi_n=R_n+S_n+o_{L^2}(1),\quad
\|R_n\|_2^2\to\alpha,\quad
\|S_n\|_2^2\to m-\alpha,
\]
for some $\alpha\in(0,m)$, with uniformly positive separation between
$\operatorname{supp}R_n$ and $\operatorname{supp}S_n$. The $L^p$ and gradient
splittings then follow exactly as in the computations below, by applying
Lemma~\ref{lem-3.2} to $(\Psi_n,R_n+S_n)$ and using the disjoint supports of
$R_n$ and $S_n$. This yields alternative \textnormal{(ii)}.

	\medskip \noindent\textbf{Step 3: $\tau=m$.}
Fix $\varepsilon>0$. By the definition of $\tau$ there exists $R>0$ such that
\[
\liminf_{n\to\infty}\rho_n(R)\ge m-\varepsilon.
\]
Passing to a subsequence if necessary, we may assume that $\rho_n(R)\ge m-\varepsilon$ for all $n$.

	For each $n\ge n_{0}$ choose $x_{n}\in\mathcal G$ with
	\[
	\int_{B_{R}(x_{n})}|\Psi_{n}(y)|^{2}\,dy
	\ge m-\varepsilon.
	\]
	By periodicity, for each $x_{n}$ one can pick $a_{n}\in\mathbb Z^{d}$ such that
	$T^{-a_{n}}x_{n}\in\mathcal K$, and define
	\[
	V_{n}=T^{-a_{n}}\Psi_{n}.
	\]
	The action is isometric, hence
	\[
	\|V_{n}\|_{H^{1}(\mathcal G)} = \|\Psi_{n}\|_{H^{1}(\mathcal G)},
	\]
	so $\{V_{n}\}$ is bounded in $H^{1}(\mathcal G)$, and
	\[
	\int_{B_{R}(T^{-a_{n}}x_{n})}|V_{n}(y)|^{2}\,dy = \int_{B_{R}(x_{n})}|\Psi_{n}(y)|^{2}\,dy
	\ge m-\varepsilon.
	\]
	Since $T^{-a_{n}}x_{n}\in\mathcal K$ and $\mathcal K$ is compact, there
	exist $x_{0}\in\mathcal K$ and a subsequence 
	such that $T^{-a_{n}}x_{n}\to x_{0}$. Choose $R'>R$ large enough that
	$B_{R}(T^{-a_{n}}x_{n})\subset B_{R'}(x_{0})$ for all sufficiently large $n$.
	Then
	\[
	\int_{B_{R'}(x_{0})}|V_{n}(y)|^{2}\,dy
	\ge m-2\varepsilon
	\]
	for all large $n$.
	
	By the boundedness of $\{V_{n}\}$ in $H^{1}(\mathcal G)$ there exist
	$\Psi\in H^{1}(\mathcal G)$ and a subsequence such that
	\[
	V_{n}\rightharpoonup\Psi\ \text{ in }H^{1}(\mathcal G),
	\]
	hence
	\[
	V_{n}\to\Psi\ \text{ in }L^{2}_{\mathrm{loc}}(\mathcal G).
	\]
	Passing to the limit in the last inequality gives
	\[
	\int_{B_{R'}(x_{0})}|\Psi(y)|^{2}\,dy
	\ge m-2\varepsilon.
	\]
	Since $\varepsilon>0$ is arbitrary,
	\[
	\|\Psi\|_{L^{2}(\mathcal G)}^{2}\ge m.
	\]
	On the other hand, weak convergence in $L^{2}(\mathcal G)$ and the lower
	semicontinuity of the norm yield
	\[
	\|\Psi\|_{L^{2}(\mathcal G)}^{2}
	\le\liminf_{n\to\infty}\|V_{n}\|_{L^{2}(\mathcal G)}^{2}
	= \lim_{n\to\infty}\|\Psi_{n}\|_{L^{2}(\mathcal G)}^{2}
	= m.
	\]
	Thus $\|\Psi\|_{L^{2}(\mathcal G)}^{2}=m$ and alternative (iii) holds.
	
	The three values $\tau=0$, $0<\tau<m$, and $\tau=m$ give alternatives (i),
	(ii), and (iii), respectively. These alternatives are mutually exclusive and
	exhaust all possibilities. This completes the proof.
\end{proof}

\section{Proof of Theorem \ref{theo-1.1}}

In order to apply the abstract critical point theorems of Section~3, we work in
the splitting $Y = Y^{-}\oplus Y^{+}$ introduced above and set
\[
M = Y^{-},\qquad N = Y^{+},
\]
so that each $u\in Y$ can be written uniquely as $u = u^{-} + u^{+}$ with
$u^{\pm}\in Y^{\pm}$. In the global, periodic setting considered here, the
action functional is
\[
\Phi(u)
= \frac{1}{2}\int_{\mathcal G}\langle u,(\mathcal D+\omega)u\rangle\,dx
- \int_{\mathcal G} F(x,u)\,dx ,
\]
so we set
\[
\Psi(u) = \int_{\mathcal G}F(x,u)\,dx .
\]

\begin{lemma}\label{lem-4.1}
	The functional $\Psi:Y\to\mathbb R$ given by
	\[
	\Psi(u) = \int_{\mathcal G}F(x,u)\,dx
	\]
	is weakly sequentially lower semicontinuous and $\Phi':Y\to Y^{*}$ is weakly
	sequentially continuous. Moreover, for every $c>0$ there exists
	$\zeta=\zeta(c)>0$ such that
	\begin{equation}\label{eq-4.1}
		\|u\| < \zeta\,\|u^{+}\|
		\qquad\text{for all }u\in\Phi_{c}.
	\end{equation}
\end{lemma}

\begin{proof}
	Let $(u_{n})\subset Y$ and $u\in Y$ be such that $u_{n}\rightharpoonup u$ in $Y$. By the continuous embeddings
	\[
	Y\hookrightarrow L^{q}(\mathcal G,\mathbb C^{2})
	\qquad\text{for all }q\in[2,\infty)
	\]
	we have $u_{n}\rightharpoonup u$ in $L^{q}(\mathcal G,\mathbb C^{2})$ for every
	such $q$. Up to a subsequence we may also assume that $u_{n}(x)\to u(x)$ for
	almost every $x\in\mathcal G$.
	
	We first prove the weak sequential lower semicontinuity of $\Psi$. We may pass to a subsequence (still denoted $(u_n)$) such that
\[
\Psi(u_n)\to \liminf_{k\to\infty}\Psi(u_k).
\] By
	$(F_{0})$–$(F_{3})$ there exist $q\in(2,\infty)$ and a constant $C_{0}>0$ such
	that
	\begin{equation}\label{eq-4.2}
		0\le F(x,z)\le C_{0}\bigl(|z|^{2}+|z|^{q}\bigr)
		\qquad\text{for all }(x,z)\in\mathcal G\times\mathbb C^{2}.
	\end{equation}
	In particular, $(F(\cdot,u_{n}(\cdot)))$ is bounded in $L^{1}(\mathcal G)$.
	
	By the $\mathbb Z^{d}$–periodicity of $\mathcal G$ there is a fundamental cell
	$\mathcal K$ such that
	\[
	\mathcal G = \bigcup_{a\in\mathbb Z^{d}}T^{a}(\mathcal K),
	\]
	with overlaps of measure zero. For $R>0$ let $\Lambda_{R}\subset\mathbb Z^{d}$
	be a finite set such that
	\[
	\mathcal G_{R}
	= \bigcup_{a\in\Lambda_{R}}T^{a}(\mathcal K)
	\]
	is connected and contains a metric ball of radius $R$. Since $\Lambda_{R}$ is
	finite and the graph automorphisms $T^{a}$ preserve the $Y$–norm, the compact
	embedding
	\[
	Y\hookrightarrow L^{q}(\mathcal K,\mathbb C^{2})
	\]
	together with a diagonal argument implies that, up to a subsequence, we have
	\begin{equation}\label{eq-4.3}
		u_{n}\to u\quad\text{in }L^{q}\bigl(T^{a}(\mathcal K),\mathbb C^{2}\bigr)
		\ \text{for every }a\in\Lambda_{R},
	\end{equation}
	and hence
	\[
	u_{n}\to u\quad\text{in }L^{q}\bigl(\mathcal G_{R},\mathbb C^{2}\bigr).
	\]
	In particular, possibly extracting a further subsequence, we may assume
	$u_{n}\to u$ almost everywhere on $\mathcal G_{R}$.
	
	Using $(F_{0})$–$(F_{3})$ again, there exists $C_{1}>0$ such that
	\[
	|F(x,u_{n}(x))|
	\le C_{1}\bigl(|u_{n}(x)|^{2}+|u_{n}(x)|^{q}\bigr)
	\qquad\text{for a.e. }x\in\mathcal G_{R}
	\]
and $C_{2}>0$ such that
	\begin{equation}\label{eq-4.4}
		\bigl|F_{u}(x,z)\bigr|
		\le C_{2}\bigl(|z|+|z|^{q-1}\bigr)
		\qquad\text{for all }(x,z)\in\mathcal G\times\mathbb C^{2}.
	\end{equation}
	Since $\mathcal G_{R}$ has finite measure and $(u_{n})$ is bounded in
	$L^{q}(\mathcal G,\mathbb C^{2})$, the right-hand side is bounded in
	$L^{1}(\mathcal G_{R})$ uniformly in $n$. Since $u_n\to u$ in $L^{q}(\mathcal G_{R},\mathbb C^{2})$ and $\mathcal G_{R}$ has finite measure,
we also have $u_n\to u$ in $L^{2}(\mathcal G_{R},\mathbb C^{2})$.
Moreover, by the mean value theorem,
\[
F(x,u_n)-F(x,u)=\int_{0}^{1}\bigl\langle F_{u}\bigl(x,u+t(u_n-u)\bigr),\,u_n-u\bigr\rangle\,dt.
\]
Using \eqref{eq-4.4} and Hölder's inequality, it follows that
\[
\|F(\cdot,u_n)-F(\cdot,u)\|_{L^{1}(\mathcal G_{R})}\to0,
\]
hence
\[
\int_{\mathcal G_{R}}F(x,u_{n})\,dx\to
\int_{\mathcal G_{R}}F(x,u)\,dx.
\]

	On the complement $\mathcal G\setminus\mathcal G_{R}$ we only use the
	nonnegativity of $F$. By Fatou's lemma,
	\[
	\int_{\mathcal G\setminus\mathcal G_{R}}F(x,u)\,dx
	\le\liminf_{n\to\infty}
	\int_{\mathcal G\setminus\mathcal G_{R}}F(x,u_{n})\,dx.
	\]
	Therefore
	\[
	\Psi(u)
	= \int_{\mathcal G}F(x,u)\,dx
	\le\liminf_{n\to\infty}\int_{\mathcal G}F(x,u_{n})\,dx
	= \liminf_{n\to\infty}\Psi(u_{n}),
	\]
	which proves the weak sequential lower semicontinuity of $\Psi$ on $Y$.
	
	We next show that $\Phi'$ is weakly sequentially continuous.
Let $(u_n)\subset Y$ be any sequence such that $u_n\rightharpoonup u$ in $Y$.
For every fixed $v\in Y$ we have
	\[
	\Phi'(u_{n})[v]
	= \bigl(\lvert\mathcal D\rvert^{1/2}u_{n}^{+},
	\lvert\mathcal D\rvert^{1/2}v^{+}\bigr)_{L^{2}}
	-\bigl(\lvert\mathcal D\rvert^{1/2}u_{n}^{-},
	\lvert\mathcal D\rvert^{1/2}v^{-}\bigr)_{L^{2}}
	+\omega\int_{\mathcal G}\langle u_{n},v\rangle\,dx
	-\int_{\mathcal G}\langle F_{u}(x,u_{n}),v\rangle\,dx.
	\]
	The linear part converges by weak convergence: since
	$u_{n}^{\pm}\rightharpoonup u^{\pm}$ in $Y^{\pm}$ and
	$\lvert\mathcal D\rvert^{1/2}:Y^{\pm}\to L^{2}(\mathcal G,\mathbb C^{2})$ is
	bounded, we obtain
	\[
	\bigl(\lvert\mathcal D\rvert^{1/2}u_{n}^{\pm},
	\lvert\mathcal D\rvert^{1/2}v^{\pm}\bigr)_{L^{2}}
	\to
	\bigl(\lvert\mathcal D\rvert^{1/2}u^{\pm},
	\lvert\mathcal D\rvert^{1/2}v^{\pm}\bigr)_{L^{2}},
	\]
	and similarly
	\[
	\int_{\mathcal G}\langle u_{n},v\rangle\,dx
	\to
	\int_{\mathcal G}\langle u,v\rangle\,dx.
	\]
	
	Thus it remains to show that
	\[
	\int_{\mathcal G}\langle F_{u}(x,u_{n}),v\rangle\,dx
	\to
	\int_{\mathcal G}\langle F_{u}(x,u),v\rangle\,dx
	\qquad\text{for every }v\in Y.
	\]
	We keep the subsets $\mathcal G_{R}$ as above. For each fixed $R>0$, the
	argument leading to \eqref{eq-4.3} gives, up to extraction of a
	subsequence, the strong convergence
	\[
	u_{n}\to u\quad\text{in }L^{q}\bigl(\mathcal G_{R},\mathbb C^{2}\bigr),
	\]
	and therefore also $u_{n}\to u$ almost everywhere on $\mathcal G_{R}$.
	
	Now we use the structural assumption $(F_{5})$. For any fixed
	$x\in\mathcal G_{R}$ and any $z_{1},z_{2}\in\mathbb C^{2}$, the integral form
	of the mean value theorem yields
	\[
	F_{u}(x,z_{1})-F_{u}(x,z_{2})
	=\int_{0}^{1}F_{uu}\bigl(x,z_{2}+t(z_{1}-z_{2})\bigr)\,(z_{1}-z_{2})\,dt.
	\]
	By $(F_{5})$ we obtain
	\[
	\bigl|F_{u}(x,z_{1})-F_{u}(x,z_{2})\bigr|
	\le c_{1}\Bigl(1+\bigl|z_{1}\bigr|^{\nu}+\bigl|z_{2}\bigr|^{\nu}\Bigr)\,\bigl|z_{1}-z_{2}\bigr|.
	\]
	Applying this with $z_{1}=u_{n}(x)$ and $z_{2}=u(x)$ gives
	\[
	\bigl|F_{u}(x,u_{n}(x))-F_{u}(x,u(x))\bigr|
	\le c_{1}\Bigl(1+|u_{n}(x)|^{\nu}+|u(x)|^{\nu}\Bigr)\,|u_{n}(x)-u(x)|.
	\]
	
Since $\mathcal G_{R}$ has finite measure and $u_n\to u$ in $L^{q}(\mathcal G_{R},\mathbb C^{2})$,
we have
\[
\|u_n-u\|_{L^{1}(\mathcal G_{R})}\le |\mathcal G_{R}|^{1-\frac{1}{q}}\|u_n-u\|_{L^{q}(\mathcal G_{R})}\to0.
\]
Moreover, for $\nu\in[0,1)$ we have $\nu\frac{q}{q-1}\le q$ (since $q>2$), hence
$|u_n|^{\nu}$ and $|u|^{\nu}$ are bounded in $L^{\frac{q}{q-1}}(\mathcal G_{R})$.
Therefore, by Hölder's inequality with exponents $\frac{q}{q-1}$ and $q$,
\[
\||u_{n}|^{\nu}(u_{n}-u)\|_{L^{1}(\mathcal G_{R})}
\le \||u_n|^{\nu}\|_{L^{\frac{q}{q-1}}(\mathcal G_{R})}\,\|u_n-u\|_{L^{q}(\mathcal G_{R})}\to0,
\]
and similarly $\||u|^{\nu}(u_n-u)\|_{L^{1}(\mathcal G_{R})}\to0$.
Consequently,
\[
\bigl\|F_{u}(\cdot,u_{n}(\cdot))-F_{u}(\cdot,u(\cdot))\bigr\|_{L^{1}(\mathcal G_{R})}\to0.
\]

	On the other hand, by \eqref{eq-4.4} and the $L^{q}$–boundedness of $(u_{n})$,
	there exists $C_{3}>0$ such that
	\[
	\|F_{u}(\cdot,u_{n}(\cdot))\|_{L^{q/(q-1)}(\mathcal G_{R})}
	\le C_{3}
	\qquad\text{for all }n,
	\]
	and the same bound holds for $F_{u}(\cdot,u(\cdot))$.
	
	Fix $v\in Y$ and define the truncation
\[
v_k(x)=
\begin{cases}
v(x), & |v(x)|\le k,\\[4pt]
k\,\frac{v(x)}{|v(x)|}, & |v(x)|>k,
\end{cases}
\]
so that $v_k\in L^{\infty}(\mathcal G_{R},\mathbb C^{2})$ and $v_k\to v$ in $L^{q}(\mathcal G_{R},\mathbb C^{2})$ as $k\to\infty$.
Then, for
	every $k$ and all $n$,
	\[
	\begin{aligned}
	\Bigl|\int_{\mathcal G_{R}}\langle F_{u}(x,u_{n}),v\rangle\,dx
	      -\int_{\mathcal G_{R}}\langle F_{u}(x,u),v\rangle\,dx\Bigr|
	&\le \Bigl|\int_{\mathcal G_{R}}\langle F_{u}(x,u_{n})-F_{u}(x,u),v_{k}\rangle\,dx\Bigr|\\
	&\quad +\Bigl|\int_{\mathcal G_{R}}\langle F_{u}(x,u_{n})-F_{u}(x,u),v-v_{k}\rangle\,dx\Bigr|.
	\end{aligned}
	\]
	The first term tends to $0$ as $n\to\infty$ for every fixed $k$, because
	$v_{k}\in L^{\infty}(\mathcal G_{R})$ and
	$F_{u}(\cdot,u_{n})\to F_{u}(\cdot,u)$ in $L^{1}(\mathcal G_{R})$. For the
	second term we use Hölder's inequality with exponents $q/(q-1)$ and $q$:
	\[
	\Bigl|\int_{\mathcal G_{R}}\langle F_{u}(x,u_{n})-F_{u}(x,u),v-v_{k}\rangle\,dx\Bigr|
	\le C_{4}\,\|v-v_{k}\|_{L^{q}(\mathcal G_{R})},
	\]
	where $C_{4}>0$ is independent of $n$ and $k$ thanks to the uniform
	$L^{q/(q-1)}$–bound on $F_{u}(\cdot,u_{n})$ and $F_{u}(\cdot,u)$. Letting
	first $n\to\infty$ and then $k\to\infty$ yields
	\[
	\int_{\mathcal G_{R}}\bigl\langle F_{u}(x,u_{n}),v\bigr\rangle\,dx
	\to
	\int_{\mathcal G_{R}}\bigl\langle F_{u}(x,u),v\bigr\rangle\,dx.
	\]
	
	On the complement $\mathcal G\setminus\mathcal G_{R}$ we estimate using
	\eqref{eq-4.4} and Hölder's inequality. Since
	$v\in Y\hookrightarrow L^{q}(\mathcal G,\mathbb C^{2})\cap L^{2}(\mathcal G,\mathbb C^{2})$,
	there exist constants $C_{5},C_{6}>0$ such that for all $n$,
	\begin{align*}
		\int_{\mathcal G\setminus\mathcal G_{R}}
		\bigl|F_{u}(x,u_{n})\cdot v\bigr|\,dx
		&\le C_{2}\int_{\mathcal G\setminus\mathcal G_{R}}
		\bigl(|u_{n}|+|u_{n}|^{q-1}\bigr)\,|v|\,dx\\
		&\le C_{5}\|u_{n}\|_{2}\,\|v\|_{L^{2}(\mathcal G\setminus\mathcal G_{R})}
		+ C_{6}\|u_{n}\|_{q}^{q-1}\,\|v\|_{L^{q}(\mathcal G\setminus\mathcal G_{R})}.
	\end{align*}
	The $Y$–boundedness of $(u_{n})$ implies that the factors $\|u_{n}\|_{2}$ and
	$\|u_{n}\|_{q}$ are uniformly bounded in $n$. Therefore, for every $\varepsilon>0$
	we can choose $R>0$ so large that
	\[
	\|v\|_{L^{2}(\mathcal G\setminus\mathcal G_{R})}
	+\|v\|_{L^{q}(\mathcal G\setminus\mathcal G_{R})} < \varepsilon,
	\]
	which yields
	\[
	\sup_{n}\int_{\mathcal G\setminus\mathcal G_{R}}
	\bigl|F_{u}(x,u_{n})\cdot v\bigr|\,dx
	\le C\,\varepsilon
	\]
	for some constant $C>0$ independent of $n$ and $R$. The same bound holds with
	$u_{n}$ replaced by $u$, since $u\in Y$ and $v$ is fixed. Hence
	\[
	\lim_{R\to\infty}\ \sup_{n}
	\int_{\mathcal G\setminus\mathcal G_{R}}
	\bigl|F_{u}(x,u_{n})\cdot v\bigr|\,dx
	= 0,
	\qquad
	\lim_{R\to\infty}
	\int_{\mathcal G\setminus\mathcal G_{R}}
	\bigl|F_{u}(x,u)\cdot v\bigr|\,dx
	= 0.
	\]
	
	Combining the convergence on $\mathcal G_{R}$ and the uniform smallness of the
	tails on $\mathcal G\setminus\mathcal G_{R}$, and letting first $n\to\infty$ and
	then $R\to\infty$, we obtain
	\[
	\int_{\mathcal G}\langle F_{u}(x,u_{n}),v\rangle\,dx
	\to
	\int_{\mathcal G}\langle F_{u}(x,u),v\rangle\,dx.
	\]
	This shows that $\Phi'(u_{n})[v]\to\Phi'(u)[v]$ for every $v\in Y$, i.e.,
	$\Phi'(u_{n})\rightharpoonup\Phi'(u)$ in $Y^{*}$. Thus $\Phi'$ is weakly
	sequentially continuous.
	
	It remains to prove the estimate \eqref{eq-4.1}. Fix $c>0$ and let
	$u\in\Phi_{c}$, so that $\Phi(u)\ge c$. Using $F\ge0$ and
	Lemma~\ref{lem-3.1}, we compute
	\begin{align}\label{eq-4.5}
		c
		&\le \Phi(u)
		= \frac{1}{2}\bigl(\|u^{+}\|^{2}-\|u^{-}\|^{2}\bigr)
		+\frac{\omega}{2}\|u\|_{2}^{2}
		-\Psi(u)\notag\\
		&\le \frac{1}{2}\bigl(\|u^{+}\|^{2}-\|u^{-}\|^{2}\bigr)
		+\frac{|\omega|}{2}\|u\|_{2}^{2}\notag\\
		&\le \frac{1}{2}\bigl(\|u^{+}\|^{2}-\|u^{-}\|^{2}\bigr)
		+\frac{|\omega|}{2m c^{2}}\|u\|^{2}\notag\\
		&= \frac{m c^{2}+|\omega|}{2m c^{2}}\|u^{+}\|^{2}
		-\frac{m c^{2}-|\omega|}{2m c^{2}}\|u^{-}\|^{2}.
	\end{align}
	Let
	\[
	\beta = \frac{m c^{2}-|\omega|}{2m c^{2}}\in(0,1),
	\qquad
	1-\beta = \frac{m c^{2}+|\omega|}{2m c^{2}}.
	\]
	Then the last line of \eqref{eq-4.5} can be written as
	\[
	c\le (1-\beta)\|u^{+}\|^{2}-\beta\|u^{-}\|^{2}.
	\]
	Since $\|u\|^{2}=\|u^{+}\|^{2}+\|u^{-}\|^{2}$, we obtain
	\[
	c\le \|u^{+}\|^{2}-\beta\|u\|^{2},
	\]
	hence
	\[
	\beta\|u\|^{2}\le\|u^{+}\|^{2}-c\le\|u^{+}\|^{2}.
	\]
	Therefore
	\[
	\|u\|^{2}\le\frac{1}{\beta}\,\|u^{+}\|^{2}
	= \frac{2m c^{2}}{m c^{2}-|\omega|}\,\|u^{+}\|^{2},
	\]
	so
	\[
	\|u\|
	\le \sqrt{\frac{2m c^{2}}{m c^{2}-|\omega|}}\ \|u^{+}\|
	=: \zeta(c)\,\|u^{+}\|.
	\]
	Consequently \eqref{eq-4.1} holds with
	$\zeta(c)=\sqrt{\frac{2m c^{2}}{m c^{2}-|\omega|}}$, and the proof is complete.
\end{proof}

\begin{lemma}\label{lem-4.2}
	There exists $\rho>0$ such that
	\[
	\varsigma:= \inf\bigl\{\Phi(u): u\in Y^{+},\ \|u\|=\rho\bigr\} > 0.
	\]
\end{lemma}

\begin{proof}
	Fix $q\in(2,\infty)$ and recall that
	\[
	\Psi(u) = \int_{\mathcal G}F(x,u)\,dx.
	\]
	By $(F_{0})$, $(F_{2})$ and $(F_{3})$, a standard growth estimate yields that
	for every $\varepsilon>0$ there exists $C_{\varepsilon}>0$ such that
	\begin{equation}\label{eq-4.6}
		F(x,u) \le \varepsilon|u|^{2} + C_{\varepsilon}|u|^{q}
		\quad\text{for all }(x,u)\in\mathcal G\times\mathbb C^{2}.
	\end{equation}
	Therefore, for all $u\in Y$,
	\begin{align}\label{eq-4.7}
		\Psi(u)
		&= \int_{\mathcal G}F(x,u)\,dx \notag\\
		&\le \varepsilon\int_{\mathcal G}|u|^{2}\,dx
		+ C_{\varepsilon}\int_{\mathcal G}|u|^{q}\,dx \notag\\
		&\le C_{1}\varepsilon\|u\|^{2}
		+ C_{2}C_{\varepsilon}\|u\|^{q},
	\end{align}
	where we used \eqref{eq-3.3} and the continuous embedding
	$Y\hookrightarrow L^{q}(\mathcal G,\mathbb C^{2})$.

	Now take $u\in Y^{+}$ with $\|u\|=\rho$. Using \eqref{eq-3.3} we obtain
	\[
	\int_{\mathcal G}|u|^{2}\,dx = \|u\|_{2}^{2}
	\le \frac{1}{m c^{2}}\|u\|^{2}
	= \frac{\rho^{2}}{m c^{2}}.
	\]
	Hence, for such $u$,
	\begin{align*}
		\Phi(u)
		&= \frac{1}{2}\|u\|^{2}
		+\frac{\omega}{2}\int_{\mathcal G}|u|^{2}\,dx
		-\Psi(u) \\
		&\ge \frac{1}{2}\|u\|^{2}
		-\frac{|\omega|}{2m c^{2}}\|u\|^{2}
		- C_{1}\varepsilon\|u\|^{2}
		- C_{2}C_{\varepsilon}\|u\|^{q} \\
		&\ge \Bigl(\frac{1}{2}
		-\frac{|\omega|}{2m c^{2}}
		- C_{1}\varepsilon\Bigr)\rho^{2}
		- C_{2}C_{\varepsilon}\rho^{q}.
	\end{align*}

        Since $|\omega|<mc^{2}$, we have $\frac{1}{2}-\frac{|\omega|}{2mc^{2}}>0$.
	Choose $\varepsilon>0$ so small that
	\[
	\mu= \frac{1}{2}
	-\frac{|\omega|}{2m c^{2}}
	- C_{1}\varepsilon > 0.
	\]
	Then choose $\rho>0$ so small that
	\[
	C_{2}C_{\varepsilon}\rho^{q-2} \le \frac{\mu}{2}.
	\]
	For such $\rho$ we obtain, for all $u\in Y^{+}$ with $\|u\|=\rho$,
	\[
	\Phi(u)
	\ge \mu\rho^{2} - C_{2}C_{\varepsilon}\rho^{q}
	\ge \frac{\mu}{2}\rho^{2} > 0.
	\]
	Therefore
	\[
	\varsigma \ge \frac{\mu}{2}\rho^{2} > 0,
	\]
	as claimed.
\end{proof}

We fix a number $\gamma$ such that
\[
m c^{2} < \gamma < b - \omega,
\]
where $b$ is given by assumption $(F_{3})$.
Let $\{E_{\lambda}\}_{\lambda\in\mathbb R}$ be the spectral family of $|\mathcal D|$ and
choose a sequence $(\gamma_{n})_{n\in\mathbb N}\subset\sigma(|\mathcal D|)\cap [m c^{2},\gamma]$ such that
\[
\gamma_{0} = m c^{2} < \gamma_{1} < \gamma_{2} < \cdots \le \gamma.
\]
For each $n\in\mathbb N$ pick an element
\[
e_{n} \in \bigl(E_{\gamma_{n}} - E_{\gamma_{n-1}}\bigr)L^{2}(\mathcal G,\mathbb C^{2})
\subset Y^{+}
\quad\text{with}\quad \|e_{n}\| = 1,
\]
and set
\[
Y_{n} = \operatorname{span}\{e_{1},\dots,e_{n}\},
\qquad
E_{n} = Y^{-} \oplus Y_{n}.
\]
By construction, the restriction of $|\mathcal D|$ to $Y_{n}$ has spectrum contained
in $[m c^{2},\gamma]$, hence for all $u^{+}\in Y_{n}$
\begin{equation}\label{eq-4.8}
	m c^{2}\,\|u^{+}\|_{2}^{2}
	\le \|u^{+}\|^{2}
	\le \gamma\,\|u^{+}\|_{2}^{2}.
\end{equation}

\begin{lemma}\label{lem-4.3}
	Assume $\omega\in(-m c^{2},m c^{2})$ and $(F_{0})$–$(F_{4})$ hold. Then for every
	$n\in\mathbb N$ one has $\sup\Phi(E_{n})<\infty$. Moreover, there exists a sequence
	$R_{n}>0$ such that
	\[
	\sup\{\Phi(u): u\in E_{n},\ \|u\|\ge R_{n}\}
	< \inf\{\Phi(u): u\in Y^{+},\ \|u\|=\rho\},
	\]
	where $\rho>0$ is given by Lemma~\ref{lem-4.2}.
\end{lemma}

\begin{proof}
	Fix $n\in\mathbb N$. We first prove that
	\begin{equation}\label{eq-4.9}
		\Phi(u)\to -\infty \quad\text{as }\|u\|\to\infty,\ u\in E_{n}.
	\end{equation}
	This will imply both $\sup\Phi(E_{n})<\infty$ and the existence of $R_{n}$ with the desired
	property.

	Suppose by contradiction that there exist $M>0$ and a sequence
	$(u_{j})\subset E_{n}$ such that
	\[
	\|u_{j}\|\to\infty,
	\qquad
	\Phi(u_{j}) \ge -M \quad\text{for all }j.
	\]
	Define the normalized sequence
	\[
	v_{j} = \frac{u_{j}}{\|u_{j}\|} \in E_{n},
	\qquad \|v_{j}\| = 1.
	\]
	Write $v_{j}=v_{j}^{-}+v_{j}^{+}$ with $v_{j}^{-} \in Y^{-}$ and $v_{j}^{+} \in Y_{n}$.
	Since $Y_{n}$ is finite dimensional and $Y^{-}$ is closed in $Y$, there exists
	$v=v^{-}+v^{+}\in E_{n}$ and a subsequence such that
	\[
	v_{j} \rightharpoonup v \text{ in }Y,\qquad
	v_{j}^{-} \rightharpoonup v^{-}\text{ in }Y^{-},\qquad
	v_{j}^{+}\to v^{+}\text{ in }Y_{n}.
	\]
	In particular $\|v\|\le1$.

	Dividing $\Phi(u_{j})$ by $\|u_{j}\|^{2}$ and using the definition of $\Phi$ we obtain
	\begin{equation}\label{eq-4.10}
		\frac{\Phi(u_{j})}{\|u_{j}\|^{2}}
		= \frac{1}{2}\Bigl(\|v_{j}^{+}\|^{2} - \|v_{j}^{-}\|^{2}
		+ \omega\|v_{j}\|_{2}^{2}\Bigr)
		- \int_{\mathcal G} \frac{F(x,u_{j})}{\|u_{j}\|^{2}}\,dx
		\ge -\frac{M}{\|u_{j}\|^{2}} = o(1).
	\end{equation}
	Using $F\ge0$ we drop the integral term and obtain
	\begin{align}\label{eq-4.11}
		o(1)
		&\le \frac{1}{2}\Bigl(\|v_{j}^{+}\|^{2} - \|v_{j}^{-}\|^{2}
		+ \omega\|v_{j}\|_{2}^{2}\Bigr) \notag\\
		&= \|v_{j}^{+}\|^{2} - \frac{1}{2}\|v_{j}\|^{2}
		+ \frac{\omega}{2}\|v_{j}\|_{2}^{2} \notag\\
		&\le \|v_{j}^{+}\|^{2}
		- \frac{1}{2}\|v_{j}\|^{2}
		+ \frac{|\omega|}{2m c^{2}}\|v_{j}\|^{2}\notag\\
		&= \|v_{j}^{+}\|^{2}
		- \frac{m c^{2}-|\omega|}{2m c^{2}}\|v_{j}\|^{2},
	\end{align}
	where we used \eqref{eq-3.3} in the third line. Since $\|v_{j}\|=1$,
	\eqref{eq-4.11} yields
	\[
	\|v_{j}^{+}\|^{2} \ge \frac{m c^{2}-|\omega|}{2m c^{2}} + o(1),
	\]
	so $\|v_{j}^{+}\|$ is bounded away from zero and therefore $v^{+}\ne 0$.

	Next we use the asymptotic behaviour of $F$ from $(F_{3})$. Define
	\[
	R(x,u) = F(x,u) - \frac{b}{2}|u|^{2}.
	\]
	By $(F_{3})$,
	\[
	\frac{R(x,u)}{|u|^{2}}\to 0 \quad\text{as }|u|\to\infty
	\]
	uniformly in $x\in\mathcal G$, and there exists a constant $C>0$ such that
	\[
	|R(x,u)| \le C|u|^{2}
	\quad\text{for all }(x,u)\in\mathcal G\times\mathbb C^{2}.
	\]
	For any $u\in Y$ we can rewrite $\Phi$ as
	\begin{align}\label{eq-4.12}
		\Phi(u)
		&= \frac{1}{2}\bigl(\|u^{+}\|^{2} - \|u^{-}\|^{2}\bigr)
		+ \frac{\omega}{2}\|u\|_{2}^{2}
		- \int_{\mathcal G}F(x,u)\,dx \notag\\
		&= \frac{1}{2}\bigl(\|u^{+}\|^{2} - \|u^{-}\|^{2}\bigr)
		+ \frac{\omega-b}{2}\|u\|_{2}^{2}
		- \int_{\mathcal G}R(x,u)\,dx.
	\end{align}
	Using the orthogonality of $Y^{+}$ and $Y^{-}$ in $L^{2}$, we have
	$\|u\|_{2}^{2} = \|u^{+}\|_{2}^{2} + \|u^{-}\|_{2}^{2}$, hence
	\begin{align}\label{eq-4.13}
		\Phi(u)
		&= \frac{1}{2}\bigl(\|u^{+}\|^{2} + \omega\|u^{+}\|_{2}^{2}\bigr)
		- \frac{1}{2}\bigl(\|u^{-}\|^{2} - \omega\|u^{-}\|_{2}^{2}\bigr)
		- \frac{b}{2}\bigl(\|u^{+}\|_{2}^{2} + \|u^{-}\|_{2}^{2}\bigr)
		- \int_{\mathcal G}R(x,u)\,dx \notag\\
		&\le \frac{1}{2}\bigl(\|u^{+}\|^{2} + \omega\|u^{+}\|_{2}^{2}\bigr)
		- \frac{m c^{2}-|\omega|}{2m c^{2}}\|u^{-}\|^{2}
		- \frac{b}{2}\|u\|_{2}^{2}
		- \int_{\mathcal G}R(x,u)\,dx,
	\end{align}
	where in the last step we used \eqref{eq-3.3} on $u^{-}$.

	We now apply \eqref{eq-4.13} to $u=v$. Using \eqref{eq-4.8} on $v^{+}$
	and \eqref{eq-3.3} on $v^{-}$, we obtain
	\begin{align}\label{eq-4.14}
		&\bigl(\|v^{+}\|^{2} + \omega\|v^{+}\|_{2}^{2}\bigr)
		- \frac{m c^{2}-|\omega|}{m c^{2}}\|v^{-}\|^{2}
		- b\|v\|_{2}^{2} \notag\\
		&\qquad\le (\gamma + \omega - b)\|v^{+}\|_{2}^{2}
		- \bigl(m c^{2}-|\omega| + b\bigr)\|v^{-}\|_{2}^{2}.
	\end{align}
	Since $\gamma < b-\omega$ and $b>0$, the coefficients on the right-hand side are
	strictly negative, and $v^{+}\ne 0$; hence
	\[
	(\gamma + \omega - b)\|v^{+}\|_{2}^{2}
	- \bigl(m c^{2}-|\omega| + b\bigr)\|v^{-}\|_{2}^{2} < 0.
	\]
	Thus the left-hand side of \eqref{eq-4.14} is strictly negative:
	\begin{equation}\label{eq-4.15}
		 \|v^{+}\|^{2} + \omega\|v^{+}\|_{2}^{2}
		- \frac{m c^{2}-|\omega|}{m c^{2}}\|v^{-}\|^{2}
		- b\|v\|_{2}^{2} < 0.
	\end{equation}

	We now localize this negativity on a bounded set. Since
	$\mathcal G$ is the union of edges of finite length, we can exhaust it by an
	increasing sequence of bounded measurable subsets
	$\Omega_{k}\subset\mathcal G$ with $\Omega_{k}\subset\Omega_{k+1}$ and
	$\bigcup_{k}\Omega_{k}=\mathcal G$. Then
	\[
	\int_{\Omega_{k}}|v|^{2}\,dx \rightarrow \int_{\mathcal G}|v|^{2}\,dx
	= \|v\|_{2}^{2}.
	\]
	Using \eqref{eq-4.15} and the monotone convergence of
	$\int_{\Omega_{k}}|v|^{2}\,dx\rightarrow\|v\|_{2}^{2}$, we can choose $k_{0}$ so large that,
	with $\Omega:=\Omega_{k_{0}}$,
	\begin{equation}\label{eq-4.16}
		\|v^{+}\|^{2} + \omega\|v^{+}\|_{2}^{2}
		- \frac{m c^{2}-|\omega|}{m c^{2}}\|v^{-}\|^{2}
		- b\int_{\Omega}|v|^{2}\,dx < 0.
	\end{equation}

	Next we estimate the nonlinear remainder on $\Omega$. For $x\in\Omega$,
	\[
	\frac{R(x,u_{j})}{\|u_{j}\|^{2}}
	= \frac{R(x,u_{j})}{|u_{j}|^{2}}\;|v_{j}(x)|^{2}
	\quad\text{whenever }u_{j}(x)\ne0,
	\]
	and we set this ratio to $0$ when $u_{j}(x)=0$.
	Fix $\varepsilon>0$ and choose $A>0$ so large that
	\[
	|u|\ge A \ \to\
	\biggl|\frac{R(x,u)}{|u|^{2}}\biggr| \le \varepsilon
	\quad\text{for all }x\in\mathcal G,
	\]
	which is possible by $(F_{3})$ and the uniformity in $x$. Then
	\[
	\Omega = \{x\in\Omega: |u_{j}(x)|\ge A\}
	\cup \{x\in\Omega: |u_{j}(x)|<A\}
	=: \Omega_{j}^{1}\cup\Omega_{j}^{2}.
	\]
	On $\Omega_{j}^{1}$ we have
	\[
	\biggl|\frac{R(x,u_{j})}{\|u_{j}\|^{2}}\biggr|
	\le \varepsilon |v_{j}(x)|^{2},
	\]
	and hence
	\[
	\int_{\Omega_{j}^{1}}\biggl|\frac{R(x,u_{j})}{\|u_{j}\|^{2}}\biggr|\,dx
	\le \varepsilon\int_{\Omega}|v_{j}|^{2}\,dx
	\le \varepsilon\,\sup_{j}\|v_{j}\|_{2}^{2}
	\le C\varepsilon.
	\]
	On $\Omega_{j}^{2}$ we only use $|R(x,u)|\le C|u|^{2}$ and $|u_{j}|\le A$ to get
	\[
	\biggl|\frac{R(x,u_{j})}{\|u_{j}\|^{2}}\biggr|
	\le \frac{C A^{2}}{\|u_{j}\|^{2}},
	\]
	so
	\[
	\int_{\Omega_{j}^{2}}\biggl|\frac{R(x,u_{j})}{\|u_{j}\|^{2}}\biggr|\,dx
	\le \frac{C A^{2}|\Omega|}{\|u_{j}\|^{2}}\xrightarrow[j\to\infty]{}0.
	\]
	Therefore
	\[
	\limsup_{j\to\infty}\int_{\Omega}\biggl|\frac{R(x,u_{j})}{\|u_{j}\|^{2}}\biggr|\,dx
	\le C\varepsilon.
	\]
	Since $\varepsilon>0$ is arbitrary, we conclude
	\begin{equation}\label{eq4.17}
		\int_{\Omega}\frac{R(x,u_{j})}{\|u_{j}\|^{2}}\,dx \to 0
		\quad\text{as }j\to\infty.
	\end{equation}

We now combine these ingredients. From \eqref{eq-4.13} applied to $u_j$ and
dividing by $\|u_j\|^{2}$, we obtain
\[
\begin{aligned}
\frac{\Phi(u_{j})}{\|u_{j}\|^{2}}
&\le \frac{1}{2}\bigl(\|v_{j}^{+}\|^{2}
+ \omega\|v_{j}^{+}\|_{2}^{2}\bigr)
- \frac{m c^{2}-|\omega|}{2m c^{2}}\|v_{j}^{-}\|^{2}
- \frac{b}{2}\int_{\mathcal G}|v_{j}|^{2}\,dx
- \int_{\mathcal G}\frac{R(x,u_{j})}{\|u_{j}\|^{2}}\,dx .
\end{aligned}
\]
Split the last two integrals over $\Omega$ and $\mathcal G\setminus\Omega$. Using
$F\ge 0$ and $R=F-\frac b2|u|^{2}$, we have
\[
-\frac{b}{2}\int_{\mathcal G\setminus\Omega}|v_{j}|^{2}\,dx
-\int_{\mathcal G\setminus\Omega}\frac{R(x,u_{j})}{\|u_{j}\|^{2}}\,dx
=
-\int_{\mathcal G\setminus\Omega}\frac{F(x,u_{j})}{\|u_{j}\|^{2}}\,dx
\le 0,
\]
hence
\[
\frac{\Phi(u_{j})}{\|u_{j}\|^{2}}
\le \frac{1}{2}\bigl(\|v_{j}^{+}\|^{2}
+ \omega\|v_{j}^{+}\|_{2}^{2}\bigr)
- \frac{m c^{2}-|\omega|}{2m c^{2}}\|v_{j}^{-}\|^{2}
- \frac{b}{2}\int_{\Omega}|v_{j}|^{2}\,dx
- \int_{\Omega}\frac{R(x,u_{j})}{\|u_{j}\|^{2}}\,dx .
\]

Since $\Omega$ is bounded (a finite union of compact edges), the embedding
$Y\hookrightarrow L^{2}(\Omega,\mathbb C^{2})$ is compact. Therefore,
$v_{j}\to v$ in $L^{2}(\Omega)$, and in particular
\[
\int_{\Omega}|v_{j}|^{2}\,dx \to \int_{\Omega}|v|^{2}\,dx .
\]
Moreover, by \eqref{eq4.17} we have
\[
\int_{\Omega}\frac{R(x,u_{j})}{\|u_{j}\|^{2}}\,dx \to 0 .
\]
Passing to $\liminf_{j\to\infty}$ and using $v_{j}^{+}\to v^{+}$ in $Y_{n}$ and
$\|v^{-}\|^{2}\le \liminf_{j\to\infty}\|v_{j}^{-}\|^{2}$, we obtain
\[
\begin{aligned}
0
\le \liminf_{j\to\infty}\frac{\Phi(u_{j})}{\|u_{j}\|^{2}}
&\le \frac{1}{2}\Bigl(
\|v^{+}\|^{2} + \omega\|v^{+}\|_{2}^{2}
- \frac{m c^{2}-|\omega|}{m c^{2}}\|v^{-}\|^{2}
- b\int_{\Omega}|v|^{2}\,dx
\Bigr).
\end{aligned}
\]
The quantity in parentheses is strictly negative by \eqref{eq-4.16}, which is a
contradiction. Hence \eqref{eq-4.9} holds.

	In particular, $\sup\Phi(E_{n})<\infty$ for each fixed $n$.
	Moreover, by Lemma~\ref{lem-4.2} there exists $\rho>0$ such that
	\[
	\varsigma= \inf\{\Phi(u): u\in Y^{+},\ \|u\|=\rho\} > 0.
	\]
	By \eqref{eq-4.9}, for each $n$ we can choose $R_{n}>0$ so large that
	\[
	\sup\{\Phi(u): u\in E_{n},\ \|u\|\ge R_{n}\} < \varsigma.
	\]
	Equivalently,
	\[
	\sup\Phi(E_{n}\setminus B_{n})
	< \inf\{\Phi(u): u\in Y^{+},\ \|u\|=\rho\},
	\]
	where $B_{n} = \{u\in E_{n}:\|u\|\le R_{n}\}$. The proof is complete.
\end{proof}

As a consequence, we have the following geometric lemma.

\begin{lemma}\label{lem-4.4}
	Assume $\omega\in(-m c^{2},m c^{2})$. There exists $R_{1}>0$ such that, for
	\[
	Q=\left\{u=u^{-}+s e_{1}:\ u^{-}\in Y^{-},\ s\ge 0,\ \|u\|\le R_{1}\right\},
	\]
	one has $\Phi\le 0$ on $\partial Q$.
\end{lemma}

\begin{proof}
	By the assumptions on $F$ we have $F(x,u)\ge 0$ for all $(x,u)$, hence
	\[
	\Psi(u)=\int_{\mathcal G}F(x,u)\,dx\ge 0.
	\]
	First consider $u\in Y^{-}$. Then $u^{+}=0$ and
	\[
	\Phi(u)
	=-\frac12\|u\|^{2}+\frac{\omega}{2}\|u\|_{2}^{2}-\Psi(u).
	\]
	By Lemma~\ref{lem-3.1} we have $m c^{2}\|u\|_{2}^{2}\le\|u\|^{2}$, so
	\[
	\Phi(u)
	\le -\frac12 m c^{2}\|u\|_{2}^{2}+\frac{\omega}{2}\|u\|_{2}^{2}-\Psi(u)
	= -\frac{m c^{2}-\omega}{2}\|u\|_{2}^{2}-\Psi(u)\le 0,
	\]
	since $m c^{2}-\omega>0$ for $\omega\in(-m c^{2},m c^{2})$ and $\Psi(u)\ge 0$.
	
	Next, consider the set $Q\subset E_{1}=Y^{-}\oplus Y_{1}$, where $Y_{1}=\mathrm{span}\{e_{1}\}$.
	By definition,
	\[
	Q = \{u\in E_{1}:\ u=u^{-}+s e_{1},\ u^{-}\in Y^{-},\ s\ge 0,\ \|u\|\le R_{1}\}.
	\]
	Thus $Q$ is the intersection of the closed ball $\{u\in E_{1}:\|u\|\le R_{1}\}$ with the
	closed half-space $\{u\in E_{1}: u=u^{-}+s e_{1},\ s\ge 0\}$. Its boundary in $E_{1}$ is the
	union of:
	\[
	\Bigl\{u=u^{-}+s e_{1}:\ s\ge 0,\ \|u\|=R_{1}\Bigr\}
	\quad\text{and}\quad
	\Bigl\{u=u^{-}\in Y^{-}:\ \|u^{-}\|\le R_{1}\Bigr\},
	\]
	corresponding to the “top” ($\|u\|=R_{1}$) and the “bottom” ($s=0$) of the cylinder.
	On the second part, $u\in Y^{-}$, so we already proved $\Phi(u)\le 0$.
	
	For the first part, note that any $u\in Q$ with $\|u\|=R_{1}$ belongs to $E_{1}=Y^{-}\oplus Y_{1}$.
	By Lemma~\ref{lem-4.3} with $n=1$ we know that
	\[
	\Phi(u)\to -\infty \quad\text{as }\|u\|\to\infty,\ u\in E_{1}.
	\]
	Hence we can choose $R_{1}>0$ so large that
	\[
	\sup\{\Phi(u):u\in E_{1},\ \|u\|\ge R_{1}\}\le 0.
	\]
	In particular, for every $u\in Q$ with $\|u\|=R_{1}$ we have $\Phi(u)\le 0$.
	
	Combining the two cases, we conclude that $\Phi\le 0$ on $\partial Q$, as claimed.
\end{proof}

\begin{lemma}\label{lem-4.5}
Any $(C)_{c}$-sequence for $\Phi$ is bounded in $Y$.
\end{lemma}

\begin{proof}
Let $(u_n)\subset Y$ be a $(C)_c$-sequence, namely
\[
\Phi(u_n)\to c,
\qquad
(1+\|u_n\|)\,\|\Phi'(u_n)\|_{Y^*}\to 0.
\]
Define
\[
\|u\|_{\omega}^{2}
=\|u\|^{2}+\omega\bigl(\|u^{+}\|_{2}^{2}-\|u^{-}\|_{2}^{2}\bigr),
\qquad
\omega_{0}=\min\{mc^{2}+\omega,\;mc^{2}-\omega\}>0.
\]
By \eqref{eq-3.3} and $\omega\in(-mc^{2},mc^{2})$ we have the norm equivalence
\begin{equation}\label{eq:omega-norm-equiv-lem45}
\omega_{0}\,\|u\|_{2}^{2}\le \|u\|_{\omega}^{2},
\qquad
\frac{mc^{2}-|\omega|}{mc^{2}}\,\|u\|^{2}\le \|u\|_{\omega}^{2}
\le \frac{mc^{2}+|\omega|}{mc^{2}}\,\|u\|^{2}
\qquad\forall\,u\in Y.
\end{equation}
In particular, $\|\cdot\|_{\omega}$ is equivalent to $\|\cdot\|$ on $Y$.

We prove that $(u_n)$ is bounded in $\|\cdot\|_\omega$. Suppose by contradiction that
\[
\|u_n\|_{\omega}\to\infty.
\]
Set
\[
v_n=\frac{u_n}{\|u_n\|_{\omega}},
\qquad \|v_n\|_{\omega}=1.
\]
Then $(v_n)$ is bounded in $Y$, and moreover $\|v_n\|_2^2\le \omega_0^{-1}$ by
\eqref{eq:omega-norm-equiv-lem45}.

\textbf{Step 1.} Testing $\Phi'(u_n)$ by $u_n^{+}-u_n^{-}$ gives
\[
\Phi'(u_n)[u_n^{+}-u_n^{-}]
=\|u_n\|_{\omega}^{2}
-\int_{\mathcal G}\bigl\langle F_u(x,u_n),\,u_n^{+}-u_n^{-}\bigr\rangle\,dx.
\]
Since $\|u_n^{+}-u_n^{-}\|\le \|u_n\|$ and
\[
|\Phi'(u_n)[u_n^{+}-u_n^{-}]|
\le \|\Phi'(u_n)\|_{Y^*}\,\|u_n^{+}-u_n^{-}\|
\le \|\Phi'(u_n)\|_{Y^*}\,\|u_n\|
\le (1+\|u_n\|)\|\Phi'(u_n)\|_{Y^*}\to 0,
\]
dividing by $\|u_n\|_\omega^2$ yields
\begin{equation}\label{eq:J_n_to_1-lem45}
J_n=\int_{\mathcal G}
\left\langle \frac{F_u(x,u_n)}{\|u_n\|_{\omega}},\,v_n^{+}-v_n^{-}\right\rangle dx
=1-\frac{\Phi'(u_n)[u_n^{+}-u_n^{-}]}{\|u_n\|_\omega^2}\to 1.
\end{equation}

\textbf{Step 2.} Assume that $(v_n)$ is vanishing in the sense of the periodic concentration--compactness
lemma (as in Lemma~\ref{Lem-3.3}). Then
\begin{equation}\label{eq:vanishing-Lp-lem45}
v_n\to 0 \quad\text{in }L^{p}(\mathcal G,\mathbb C^{2})
\quad\text{for every }p\in(2,\infty).
\end{equation}
Define
\[
\hat F(x,u)=\frac{1}{2}\langle F_u(x,u),u\rangle - F(x,u).
\]
Using the identity
\[
\Phi(u)-\frac12\Phi'(u)[u]=\int_{\mathcal G}\hat F(x,u)\,dx,
\]
the Cerami condition gives
\[
\int_{\mathcal G}\hat F(x,u_n)\,dx
=\Phi(u_n)-\frac12\Phi'(u_n)[u_n]\to c,
\]
hence $\int_{\mathcal G}\hat F(x,u_n)\,dx$ is bounded.
Recall that $(F_4)$ yields $\hat F\ge 0$ and provides $\delta_1\in(0,\omega_0)$ such that
\[
\hat F(x,u)\ge\delta_1
\quad\text{whenever}\quad
|F_u(x,u)|\ge (\omega_0-\delta_1)|u|.
\]
Set
\[
S_n=\Bigl\{x\in\mathcal G:\ |F_u(x,u_n(x))|\ge (\omega_0-\delta_1)|u_n(x)|\Bigr\},
\qquad
T_n=\mathcal G\setminus S_n.
\]
Then $\hat F(x,u_n)\ge\delta_1$ on $S_n$, hence
\begin{equation}\label{eq:Sn-measure-lem45}
|S_n|\le \frac{1}{\delta_1}\int_{\mathcal G}\hat F(x,u_n)\,dx \le C
\quad\text{for all }n.
\end{equation}

We estimate $J_n$ by splitting $S_n$ and $T_n$.
On $T_n$ we have $|F_u(x,u_n)|\le (\omega_0-\delta_1)|u_n|=(\omega_0-\delta_1)\|u_n\|_\omega |v_n|$, hence
\[
\left|\int_{T_n}
\left\langle \frac{F_u(x,u_n)}{\|u_n\|_{\omega}},\,v_n^{+}-v_n^{-}\right\rangle dx\right|
\le (\omega_0-\delta_1)\int_{\mathcal G}|v_n|\,|v_n^{+}-v_n^{-}|\,dx.
\]
By Cauchy--Schwarz and the $L^2$-orthogonality of $Y^+$ and $Y^-$,
\[
\int_{\mathcal G}|v_n|\,|v_n^{+}-v_n^{-}|
\le \|v_n\|_2\,\|v_n^{+}-v_n^{-}\|_2
=\|v_n\|_2^2
\le \frac{1}{\omega_0},
\]
so
\begin{equation}\label{eq:Tn-bound-lem45}
\limsup_{n\to\infty}\left|\int_{T_n}\left\langle \frac{F_u(x,u_n)}{\|u_n\|_{\omega}},\,v_n^{+}-v_n^{-}\right\rangle dx\right|
\le \frac{\omega_0-\delta_1}{\omega_0}=1-\frac{\delta_1}{\omega_0}.
\end{equation}

On $S_n$ we use a global growth bound for $F_u$.
From $(F_3)$ and the Carath\'eodory regularity of $F_u$, there exists $C_*>0$ such that
\begin{equation}\label{eq:Fu-growth-lem45}
|F_u(x,z)|\le C_*(1+|z|)
\qquad\text{for all }(x,z)\in\mathcal G\times\mathbb C^2.
\end{equation}
Therefore, for all $n$,
\[
\left|\int_{S_n}
\left\langle \frac{F_u(x,u_n)}{\|u_n\|_{\omega}},\,v_n^{+}-v_n^{-}\right\rangle dx\right|
\le \left\|\frac{F_u(\cdot,u_n)}{\|u_n\|_\omega}\right\|_{L^2(S_n)}
\|v_n^{+}-v_n^{-}\|_{L^2(S_n)}.
\]
Using \eqref{eq:Fu-growth-lem45} and $\|v_n^{+}-v_n^{-}\|_2=\|v_n\|_2\le \omega_0^{-1/2}$,
we get
\[
\left\|\frac{F_u(\cdot,u_n)}{\|u_n\|_\omega}\right\|_{L^2(S_n)}
\le C_*\left(\|v_n\|_{L^2(S_n)}+\frac{|S_n|^{1/2}}{\|u_n\|_\omega}\right),
\]
hence
\begin{equation}\label{eq:Sn-bound-lem45}
\left|\int_{S_n}
\left\langle \frac{F_u(x,u_n)}{\|u_n\|_{\omega}},\,v_n^{+}-v_n^{-}\right\rangle dx\right|
\le C\left(\|v_n\|_{L^2(S_n)}+\frac{|S_n|^{1/2}}{\|u_n\|_\omega}\right).
\end{equation}
By \eqref{eq:Sn-measure-lem45} and \eqref{eq:vanishing-Lp-lem45}, choosing any $p>2$ we have
\[
\|v_n\|_{L^2(S_n)}
\le |S_n|^{\frac12-\frac1p}\|v_n\|_{L^p(\mathcal G)}\to 0,
\]
and also $|S_n|^{1/2}/\|u_n\|_\omega\to 0$ since $|S_n|$ is uniformly bounded and
$\|u_n\|_\omega\to\infty$. Thus the right-hand side of \eqref{eq:Sn-bound-lem45} tends to $0$, and
\[
\lim_{n\to\infty}\int_{S_n}
\left\langle \frac{F_u(x,u_n)}{\|u_n\|_{\omega}},\,v_n^{+}-v_n^{-}\right\rangle dx=0.
\]
Combining with \eqref{eq:Tn-bound-lem45} gives
\[
\limsup_{n\to\infty}J_n \le 1-\frac{\delta_1}{\omega_0}<1,
\]
which contradicts \eqref{eq:J_n_to_1-lem45}. Hence vanishing cannot occur.

\textbf{Step 3.} Since vanishing is excluded, by Lemma \ref{Lem-3.3} there exist
a sequence of shifts $T^{a_n}$ such that
\[
\tilde v_n=T^{-a_n}v_n \rightharpoonup \tilde v \text{ in }Y,
\qquad
\tilde v_n\to \tilde v \text{ in }L^2_{\mathrm{loc}}(\mathcal G,\mathbb C^2),
\qquad
\tilde v\not\equiv 0.
\]
Set $\tilde u_n=T^{-a_n}u_n=\|u_n\|_\omega \tilde v_n$.

Fix $\varphi\in C_c^\infty(\mathcal G,\mathbb C^2)$ and define $\varphi_n=T^{a_n}\varphi$.
By periodicity of $\Phi'$ and the isometric action of $T^{a}$ on $Y$,
\[
\Phi'(u_n)[\varphi_n]=\Phi'(\tilde u_n)[\varphi],
\qquad
\|\varphi_n\|=\|\varphi\|.
\]
Moreover,
\[
\left|\frac{\Phi'(\tilde u_n)[\varphi]}{\|u_n\|_\omega}\right|
=\left|\frac{\Phi'(u_n)[\varphi_n]}{\|u_n\|_\omega}\right|
\le \frac{\|\Phi'(u_n)\|_{Y^*}\,\|\varphi_n\|}{\|u_n\|_\omega}
\le \frac{\|\varphi\|}{\|u_n\|_\omega}\,\|\Phi'(u_n)\|_{Y^*}\to 0,
\]
because $\|u_n\|_\omega\to\infty$ and $\|\Phi'(u_n)\|_{Y^*}\to 0$.

Writing out $\Phi'$ and dividing by $\|u_n\|_\omega$ yields
\[
\begin{aligned}
0
&=\lim_{n\to\infty}\frac{\Phi'(\tilde u_n)[\varphi]}{\|u_n\|_\omega}\\
&=\lim_{n\to\infty}\Bigl[
(|\mathcal D|^{1/2}\tilde v_n^{+},|\mathcal D|^{1/2}\varphi^{+})_{L^2}
-(|\mathcal D|^{1/2}\tilde v_n^{-},|\mathcal D|^{1/2}\varphi^{-})_{L^2}
+\omega(\tilde v_n,\varphi)_{L^2}
-\int_{\mathcal G}\Bigl\langle \frac{F_u(x,\tilde u_n)}{\|u_n\|_\omega},\varphi\Bigr\rangle dx
\Bigr].
\end{aligned}
\]
The linear terms converge to
\[
(|\mathcal D|^{1/2}\tilde v^{+},|\mathcal D|^{1/2}\varphi^{+})_{L^2}
-(|\mathcal D|^{1/2}\tilde v^{-},|\mathcal D|^{1/2}\varphi^{-})_{L^2}
+\omega(\tilde v,\varphi)_{L^2}.
\]
It remains to identify the nonlinear limit. Let $\Omega=\mathrm{supp}\,\varphi$.
Fix $\varepsilon>0$ and choose $R>0$ so large that by $(F_3)$,
\[
|F_u(x,z)-bz|\le \varepsilon|z|
\quad\text{for all }x\in\mathcal G,\ |z|\ge R.
\]
Split $\Omega=\Omega_n^1\cup\Omega_n^0$ with
\[
\Omega_n^1=\{x\in\Omega:\ |\tilde u_n(x)|\ge R\},
\qquad
\Omega_n^0=\Omega\setminus\Omega_n^1.
\]
On $\Omega_n^1$ we write $F_u(x,\tilde u_n)=b\tilde u_n+r_n$ with $|r_n|\le \varepsilon|\tilde u_n|$, hence
\[
\frac{F_u(x,\tilde u_n)}{\|u_n\|_\omega}=b\tilde v_n+\tilde r_n,
\qquad
|\tilde r_n|\le \varepsilon|\tilde v_n|.
\]
Since $\Omega$ is bounded, the embedding $Y\hookrightarrow L^2(\Omega)$ is compact, so
$\tilde v_n\to \tilde v$ in $L^2(\Omega)$. Therefore,
\[
\int_{\Omega_n^1}\langle b\tilde v_n,\varphi\rangle dx
= b(\tilde v_n,\varphi)_{L^2(\Omega)}-b(\tilde v_n,\varphi)_{L^2(\Omega_n^0)}
\to b(\tilde v,\varphi)_{L^2(\Omega)},
\]
because on $\Omega_n^0$ we have $|\tilde v_n|\le R/\|u_n\|_\omega\to 0$ uniformly, hence
$(\tilde v_n,\varphi)_{L^2(\Omega_n^0)}\to 0$.
Moreover,
\[
\left|\int_{\Omega_n^1}\langle \tilde r_n,\varphi\rangle dx\right|
\le \varepsilon\|\tilde v_n\|_{L^2(\Omega)}\|\varphi\|_{L^2(\Omega)}
\le C\varepsilon,
\]
with $C$ independent of $n$.
On $\Omega_n^0$ we have $|\tilde u_n|\le R$, hence by continuity and periodicity there exists $M_R>0$
such that $|F_u(x,z)|\le M_R$ for all $x\in\mathcal G$, $|z|\le R$, and thus
\[
\left|\int_{\Omega_n^0}\left\langle \frac{F_u(x,\tilde u_n)}{\|u_n\|_\omega},\varphi\right\rangle dx\right|
\le \frac{M_R}{\|u_n\|_\omega}\|\varphi\|_{L^1(\Omega)}\to 0.
\]
Combining the pieces and letting $\varepsilon\to 0$ yields
\[
\int_{\mathcal G}\Bigl\langle \frac{F_u(x,\tilde u_n)}{\|u_n\|_\omega},\varphi\Bigr\rangle dx
\to b(\tilde v,\varphi)_{L^2(\mathcal G)}.
\]
Therefore,
\[
(|\mathcal D|^{1/2}\tilde v^{+},|\mathcal D|^{1/2}\varphi^{+})_{L^2}
-(|\mathcal D|^{1/2}\tilde v^{-},|\mathcal D|^{1/2}\varphi^{-})_{L^2}
+\omega(\tilde v,\varphi)_{L^2}
-b(\tilde v,\varphi)_{L^2}=0
\qquad\forall\,\varphi\in C_c^\infty(\mathcal G,\mathbb C^2).
\]
By density of $C_c^\infty(\mathcal G,\mathbb C^2)$ in $Y$ and continuity of the above identity in $\varphi$,
it holds for every $\varphi\in Y$. By the representation theorem for closed forms associated with the
self-adjoint operator $\mathcal D$, this implies $\tilde v\in\mathrm{dom}(\mathcal D)$ and
\[
\mathcal D\tilde v=(b-\omega)\tilde v \quad\text{in }L^2(\mathcal G,\mathbb C^2).
\]
Since $\tilde v\not\equiv 0$, we obtain $b-\omega\in\sigma_p(\mathcal D)$, which contradicts the
additional requirement in $(F_3)$ that $b-\omega\notin\sigma_p(\mathcal D)$.

Thus the assumption $\|u_n\|_\omega\to\infty$ is false, so $(u_n)$ is bounded in $\|\cdot\|_\omega$.
Finally, \eqref{eq:omega-norm-equiv-lem45} implies that $(u_n)$ is bounded in $\|\cdot\|$, hence bounded in $Y$.
\end{proof}

Let
\[
\mathcal{C}
=\{u\in Y\setminus\{0\}:\ \Phi'(u)=0\}
\]
be the set of nontrivial critical points of \(\Phi\).
We argue by contradiction and assume that
\begin{equation}\label{eq-4.22}
	\mathcal{C}/\mathbb{Z}^{d} \text{ is a finite set}.
\end{equation}
Here the quotient is taken with respect to the \(\mathbb{Z}^{d}\)–action induced
by the graph isometries \(\{T^{k}\}_{k\in\mathbb{Z}^{d}}\), that is,
\(u\sim v\) if and only if \(v=T^{-k}u\) for some \(k\in\mathbb{Z}^{d}\).
We will show that under \eqref{eq-4.22} condition \((\Phi_{5})\) is satisfied.
Then Theorem~\ref{theo-3.2} yields an unbounded sequence of positive critical
values of \(\Phi\), which contradicts \eqref{eq-4.22}.

Under the assumption \eqref{eq-4.22}, let \(\mathcal{F}\) be a finite set of
representatives of the \(\mathbb{Z}^{d}\)–orbits in \(\mathcal{C}\), so that each
\(u\in\mathcal{C}\) can be written as \(u=T^{-k}w\) for some \(k\in\mathbb{Z}^{d}\)
and some \(w\in\mathcal{F}\).
Since \(\Phi\) is even (by the assumptions on \(F\))
and \(\Phi'\) is odd, we may assume that \(\mathcal{F}\) is symmetric, that is
\(\mathcal{F}=-\mathcal{F}\).

If $u\in\mathcal C$ and $u\neq 0$, then $\Phi(u)>0$. Indeed, since $\Phi'(u)=0$,
\[
\Phi(u)=\Phi(u)-\frac12\Phi'(u)[u]=\int_{\mathcal G}\hat F(x,u)\,dx\ge 0,
\qquad
\hat F(x,z)=\frac12\langle F_u(x,z),z\rangle-F(x,z),
\]
and $\hat F\ge 0$ by $(F_4)$. Let $\delta_1\in(0,\omega_0)$ be given by $(F_4)$ and set
\[
S=\Bigl\{x\in\mathcal G:\ |F_u(x,u(x))|\ge (\omega_0-\delta_1)|u(x)|\Bigr\}.
\]
Then $\hat F(x,u(x))\ge \delta_1$ a.e.\ on $S$. If $|S|=0$, testing $\Phi'(u)=0$ with
$u^+-u^-$ gives $\|u\|_\omega^2=\int_{\mathcal G}\langle F_u(x,u),u^+-u^-\rangle dx$ and hence
\[
\|u\|_\omega^2\le (\omega_0-\delta_1)\|u\|_2^2,
\]
while \eqref{eq:omega-norm-equiv-lem45} yields $\|u\|_\omega^2\ge \omega_0\|u\|_2^2$, forcing $u=0$,
a contradiction. Thus $|S|>0$ and
\[
\Phi(u)=\int_{\mathcal G}\hat F(x,u)\,dx\ge \int_S\hat F(x,u)\,dx\ge \delta_1|S|>0.
\]

Since \(\mathcal{F}\) is finite, there exist constants \(0<\theta\le\vartheta\)
such that
\begin{equation}\label{eq-4.23}
	\theta
	<\min_{w\in\mathcal{F}}\Phi(w)
	=\min_{u\in\mathcal{C}}\Phi(u)
	\le\max_{u\in\mathcal{C}}\Phi(u)
	=\max_{w\in\mathcal{F}}\Phi(w)
	<\vartheta.
\end{equation}
Let \([r]\) denote the integer part of \(r\in\mathbb{R}\).

\begin{lemma}\label{lem-4.6}
	Assume \eqref{eq-4.22} holds and let $(u_{m})$ be a $(C)_{c}$-sequence for $\Phi$ in $Y$.
	Then either
	\begin{itemize}
		\item[(i)] $u_{m}\to0$ in $Y$ and $c=0$, or
		\item[(ii)] $c\ge\theta$ and there exist a positive integer
		$\ell\le[c/\theta]$, points $\bar u_{1},\ldots,\bar u_{\ell}\in\mathcal{F}$, a subsequence (still denoted by $(u_m)$), and sequences
		$(a_{m}^{i})\subset\mathbb{Z}^{d}$, $i=1,\ldots,\ell$, such that
		\[
		\Bigl\|u_{m}-\sum_{i=1}^{\ell} (a_{m}^{i}* \bar u_{i})\Bigr\|\;\to\;0
		\quad\text{and}\quad
		\sum_{i=1}^{\ell}\Phi(\bar u_{i})=c,
		\]
		where, for $a\in\mathbb{Z}^{d}$ and $u\in Y$, the translate
		\((a*u)\) is given by
		\[
		(a*u)(x)=u\bigl(T^{-a}x\bigr).
		\]
	\end{itemize}
\end{lemma}

\begin{proof}
By Lemma~\ref{lem-4.5}, $(u_m)$ is bounded in $Y$. Using $(F_4)$ and
$\Phi'(u_m)\to 0$ in $Y^*$, we have
\[
0\le \int_{\mathcal G}\hat F(x,u_m)\,dx
=\Phi(u_m)-\frac12\Phi'(u_m)[u_m]\to c,
\]
hence $c\ge 0$.

If $u_m\to 0$ in $Y$, then $\int_{\mathcal G}\hat F(x,u_m)\,dx\to 0$ and thus $c=0$,
which is (i). From now on assume that $u_m\not\to 0$ in $Y$.

Recall
\[
\|u\|_\omega^2=\|u\|^2+\omega\bigl(\|u^+\|_2^2-\|u^-\|_2^2\bigr),
\qquad
\omega_0=\min\{mc^2-\omega,mc^2+\omega\}>0,
\]
and the equivalence
\[
\omega_0\|u\|_2^2\le \|u\|_\omega^2,
\qquad
\frac{mc^2-|\omega|}{mc^2}\|u\|^2\le \|u\|_\omega^2
\le \frac{mc^2+|\omega|}{mc^2}\|u\|^2 .
\]

\textbf{Step 1.}
Assume that $(u_m)$ is vanishing in the sense of Lemma~\ref{Lem-3.3}. Then
$\|u_m\|_{L^p(\mathcal G)}\to 0$ for every $p\in(2,\infty)$.
Fix one such $p$. By $(F_2)$--$(F_3)$, for every $\varepsilon>0$ there exists
$C_\varepsilon>0$ such that
\[
|F_u(x,z)|\le \varepsilon |z|+C_\varepsilon |z|^{p-1}
\qquad\text{for all }(x,z)\in \mathcal G\times \mathbb C^2.
\]
Testing $\Phi'(u_m)$ by $u_m^+-u_m^-$ gives
\[
\Phi'(u_m)[u_m^+-u_m^-]
=\|u_m\|_\omega^2-\int_{\mathcal G}\langle F_u(x,u_m),u_m^+-u_m^-\rangle\,dx.
\]
Since $(u_m)$ is a $(C)_c$-sequence, $\Phi'(u_m)[u_m^+-u_m^-]\to 0$, hence
\[
\|u_m\|_\omega^2
=\int_{\mathcal G}\langle F_u(x,u_m),u_m^+-u_m^-\rangle\,dx+o(1).
\]
By the above growth bound, H\"older inequality, and the boundedness of $(u_m)$ in $Y$, we get
\[
\left|\int_{\mathcal G}\langle F_u(x,u_m),u_m^+-u_m^-\rangle\,dx\right|
\le \varepsilon C_1 + C_\varepsilon C_2\,\|u_m\|_{L^p(\mathcal G)}^{p-1}\to \varepsilon C_1,
\]
with constants $C_1,C_2$ independent of $m$. Letting $\varepsilon\to 0$ yields
$\|u_m\|_\omega\to 0$, hence $\|u_m\|\to 0$ by norm equivalence, contradicting
$u_m\not\to 0$ in $Y$. Therefore $(u_m)$ is nonvanishing.

\textbf{Step 2.}
By nonvanishing and Lemma~\ref{Lem-3.3}, there exist a sequence $(b_m^1)\subset \mathbb Z^d$
and $v^{(1)}\in Y$, $v^{(1)}\neq 0$, such that 
\[
b_m^1*u_m \rightharpoonup v^{(1)} \quad\text{in }Y,
\qquad
b_m^1*u_m \to v^{(1)} \quad\text{in }L^p_{\mathrm{loc}}(\mathcal G)
\ \text{for all }p\in[2,\infty).
\]
Since $\Phi$ and $\Phi'$ are $\mathbb Z^d$-invariant, $(b_m^1*u_m)$ is still a $(C)_c$-sequence.
By Lemma~\ref{lem-4.1}, we obtain $\Phi'(v^{(1)})=0$,
so $v^{(1)}\in \mathcal C$.

Choose $\bar u_1\in\mathcal F$ and $k_1\in\mathbb Z^d$ such that $v^{(1)}=k_1*\bar u_1$,
and replace $b_m^1$ with $b_m^1+k_1$ so that $v^{(1)}=\bar u_1$.
Set $a_m^1=-b_m^1$ and define the remainder
\[
r_m^{(1)}=u_m-a_m^1*\bar u_1.
\]
Then
\[
b_m^1*r_m^{(1)}=b_m^1*u_m-\bar u_1 \rightharpoonup 0\quad\text{in }Y,
\qquad
b_m^1*r_m^{(1)}\to 0\quad\text{in }L^p_{\mathrm{loc}}(\mathcal G).
\]
Using Lemma~\ref{lem-3.2} and the $\mathbb Z^d$-invariance of $\Phi$, we get
\begin{equation}\label{eq:split-1}
\Phi(u_m)=\Phi(\bar u_1)+\Phi(r_m^{(1)})+o(1),
\qquad
\Phi'(r_m^{(1)})\to 0.
\end{equation}
In particular, $\Phi(\bar u_1)\ge \theta$ by \eqref{eq-4.23}, hence $c\ge \theta$ whenever
alternative (ii) occurs.

If $c=\Phi(\bar u_1)$, then $\Phi(r_m^{(1)})\to 0$ and $\Phi'(r_m^{(1)})\to 0$.
If $r_m^{(1)}\not\to 0$ in $Y$, then $r_m^{(1)}$ is nonvanishing and Step 2 can be repeated to
extract another profile $\bar u_2\in\mathcal F$ with $\Phi(\bar u_2)\ge\theta$, which would force
$\liminf_m \Phi(r_m^{(1)})\ge \theta$, a contradiction. Hence $r_m^{(1)}\to 0$ in $Y$ and (ii)
holds with $\ell=1$.

\textbf{Step 3.}
Assume $c>\Phi(\bar u_1)$ and set $c_1=c-\Phi(\bar u_1)>0$. Then by \eqref{eq:split-1},
$(r_m^{(1)})$ is a bounded $(C)_{c_1}$-sequence.
If $r_m^{(1)}$ were vanishing, Step 1 would give $r_m^{(1)}\to 0$ in $Y$, hence $c_1=0$,
contradiction. Therefore $r_m^{(1)}$ is nonvanishing.

Applying Lemma~\ref{Lem-3.3} to $r_m^{(1)}$, we find $(b_m^2)\subset\mathbb Z^d$ and
$v^{(2)}\in\mathcal C\setminus\{0\}$ such that $b_m^2*r_m^{(1)}\rightharpoonup v^{(2)}$.
Moreover, because $b_m^1*r_m^{(1)}\to 0$ in $L^2_{\mathrm{loc}}$, necessarily
\[
|b_m^2-b_m^1|\to\infty,
\]
otherwise $b_m^2*r_m^{(1)}$ would also converge to $0$ locally and could not have a nontrivial limit.
Choose $\bar u_2\in\mathcal F$ as the orbit representative of $v^{(2)}$ and adjust $b_m^2$ by a fixed
shift so that $v^{(2)}=\bar u_2$. Set $a_m^2=-b_m^2$ and define
\[
r_m^{(2)}=r_m^{(1)}-a_m^2*\bar u_2.
\]
As before, Lemma~\ref{lem-3.2} and invariance yield
\[
\Phi(r_m^{(1)})=\Phi(\bar u_2)+\Phi(r_m^{(2)})+o(1),
\qquad
\Phi'(r_m^{(2)})\to 0.
\]
Consequently,
\[
\Phi(u_m)=\Phi(\bar u_1)+\Phi(\bar u_2)+\Phi(r_m^{(2)})+o(1),
\qquad
c_2:=c-\Phi(\bar u_1)-\Phi(\bar u_2)\ge 0,
\]
and $|a_m^2-a_m^1|=|b_m^2-b_m^1|\to\infty$.

Iterating, we obtain profiles $\bar u_1,\dots,\bar u_k\in\mathcal F$, shifts
$(a_m^i)\subset\mathbb Z^d$ with $|a_m^i-a_m^j|\to\infty$ for $i\neq j$, and remainders $r_m^{(k)}$
such that
\[
u_m=\sum_{i=1}^k (a_m^i*\bar u_i)+r_m^{(k)},
\qquad
\Phi'(r_m^{(k)})\to 0,
\qquad
\Phi(r_m^{(k)})\to c_k:=c-\sum_{i=1}^k \Phi(\bar u_i)\ge 0.
\]
Since $\Phi(\bar u_i)\ge \theta$ for all $i$, we must have $k\le [c/\theta]$.

Let $\ell$ be the maximal number of extracted profiles, so that $c_\ell\in[0,\theta)$.
If $c_\ell>0$, then $(r_m^{(\ell)})$ is a bounded nonvanishing $(C)_{c_\ell}$-sequence and the
above procedure would produce an additional profile $\bar u_{\ell+1}\in\mathcal F$ with
$\Phi(\bar u_{\ell+1})\ge \theta$, giving
\[
c=\sum_{i=1}^{\ell+1}\Phi(\bar u_i)+c_{\ell+1}\ge \sum_{i=1}^{\ell+1}\Phi(\bar u_i)
\ge c-\theta+\theta=c,
\]
which forces $c_{\ell+1}=0$ and contradicts the maximality of $\ell$. Hence $c_\ell=0$.

Finally, if $r_m^{(\ell)}\not\to 0$ in $Y$, then $r_m^{(\ell)}$ is nonvanishing and we could extract
one more profile with energy at least $\theta$, contradicting $c_\ell=0$. Therefore
$r_m^{(\ell)}\to 0$ in $Y$, and we conclude
\[
\left\|u_m-\sum_{i=1}^\ell (a_m^i*\bar u_i)\right\|\to 0,
\qquad
\sum_{i=1}^\ell \Phi(\bar u_i)=c,
\qquad
\ell\le [c/\theta].
\]
This is alternative (ii). The proof is complete.
\end{proof}

For $\ell\in\mathbb{N}$ and a finite set $\mathcal{B} \subset Y$ we define
\[
[\mathcal{B}, \ell]
=\Bigl\{\sum_{i=1}^{j}(k_i * u_i):
\ 1 \le j \le \ell,\ k_i \in \mathbb{Z}^{d},\ u_i \in \mathcal{B}\Bigr\},
\]
where $(k*u)(x)=u\bigl(T^{-k}x\bigr)$ denotes the $\mathbb{Z}^{d}$–action induced by
the graph automorphisms $T^{k}$ on $\mathcal{G}$.

An argument similar to the one in \cite{MR1070929} shows that
\begin{equation}\label{eq-4.27}
	\inf\Bigl\{\|u-u'\|:\ u,u'\in[\mathcal{B},\ell],\ u\neq u'\Bigr\}>0.
\end{equation}

As a consequence of Lemma~\ref{lem-4.6} we have the following.

\begin{lemma}\label{lem-4.7}
	Assume \eqref{eq-4.22} holds. Then $\Phi$ satisfies $(\Phi_{5})$.
\end{lemma}

\begin{proof}
Let $I\subset(0,\infty)$ be a compact interval and set $c_\ast=\max I$.
Choose
\[
\ell=\bigl[c_\ast/\theta\bigr],
\qquad
\mathcal A=[\mathcal F,\ell],
\]
where $\mathcal F$ is the finite set of $\mathbb Z^d$--orbit representatives of $\mathcal C$.
Let $P^+:Y\to Y^+$ be the orthogonal projection.

Since the translations $T^k$ commute with $\mathcal D$, they preserve the spectral splitting
$Y=Y^-\oplus Y^+$. In particular,
\[
P^+(k*u)=k*(P^+u)\quad\text{for all }k\in\mathbb Z^d,\ u\in Y,
\]
hence
\[
P^+\mathcal A=[P^+\mathcal F,\ell].
\]
Applying \eqref{eq-4.27} with $\mathcal B=P^+\mathcal F$ yields
\[
\inf\Bigl\{\|w-w'\|:\ w,w'\in P^+\mathcal A,\ w\neq w'\Bigr\}>0,
\]
so $P^+\mathcal A$ is uniformly separated in $Y^+$.

Moreover, since each translation is an isometry on $Y$, every $u\in\mathcal A$ can be written as
\[
u=\sum_{i=1}^j (k_i*u_i),\qquad 1\le j\le \ell,\ k_i\in\mathbb Z^d,\ u_i\in\mathcal F,
\]
and therefore
\[
\|u\|\le \sum_{i=1}^j\|u_i\|
\le \ell\,\max\bigl\{\|\bar u\|:\ \bar u\in\mathcal F\bigr\}
\qquad\text{for all }u\in\mathcal A.
\]
Hence $\mathcal A$ is bounded in $Y$.
In addition, \eqref{eq-4.27} implies that $\mathcal A$ is uniformly separated in $Y$,
so $\mathcal A$ is a closed subset of $Y$.

Now fix $c\in I$ and let $(u_m)$ be a $(C)_c$--sequence for $\Phi$.
Since $c>0$, alternative (i) in Lemma~\ref{lem-4.6} cannot occur. Thus, after passing to a subsequence,
there exist an integer $\ell_c\le [c/\theta]\le \ell$, elements $\bar u_1,\ldots,\bar u_{\ell_c}\in\mathcal F$,
and sequences $(a_m^i)\subset\mathbb Z^d$ such that
\[
\Bigl\|u_m-\sum_{i=1}^{\ell_c}(a_m^i*\bar u_i)\Bigr\|\to 0,
\qquad
\sum_{i=1}^{\ell_c}\Phi(\bar u_i)=c.
\]
Set
\[
w_m:=\sum_{i=1}^{\ell_c}(a_m^i*\bar u_i).
\]
Then $w_m\in[\mathcal F,\ell_c]\subset[\mathcal F,\ell]=\mathcal A$ for all $m$, hence
\[
\operatorname{dist}(u_m,\mathcal A)\le \|u_m-w_m\|\to 0
\]
along that subsequence. Therefore, $\mathcal A$ is a bounded $(C)_I$--attractor for $\Phi$.

By definition of $(\Phi_5)$, the existence of such a bounded $(C)_I$--attractor $\mathcal A$
for every compact $I\subset(0,\infty)$, together with the uniform separation of $P^+\mathcal A$
in $Y^+$, implies that $\Phi$ satisfies $(\Phi_5)$.
\end{proof}

{\bf Proof of Theorem~\ref{theo-1.1}:}
Let $M=Y^{-}$ and $N=Y^{+}$. By Lemma~\ref{lem-4.1} and
Theorem~\ref{theo-3.3}, the functional $\Phi$ satisfies $(\Phi_{0})$ and
$(\Phi_{1})$. Lemma~\ref{lem-4.2} yields $(\Phi_{2})$, while
Lemma~\ref{lem-4.3} and Lemma~\ref{lem-4.4} provide the linking geometry
required in Theorem~\ref{theo-3.1}. Hence all the assumptions of
Theorem~\ref{theo-3.1} are fulfilled.

Therefore there exists a sequence $(u_{m})\subset Y$ such that
\[
\Phi(u_{m})\to c\ge\varsigma
\quad\text{and}\quad
\bigl(1+\|u_{m}\|\bigr)\,\|\Phi'(u_{m})\|_{Y^{*}}\to 0,
\]
that is, $(u_{m})$ is a $(C)_{c}$–sequence at some level $c\ge\varsigma>0$.
By Lemma~\ref{lem-4.5}, the sequence $(u_{m})$ is bounded in $Y$, and hence
$\Phi'(u_{m})\to 0$ in $Y^{*}$.

\textbf{Step 1.}
Assume by contradiction that $(u_m)$ is vanishing in the sense of Lemma~\ref{Lem-3.3}, i.e.
for some $r>0$,
\[
\lim_{m\to\infty}\ \sup_{x\in\mathcal G}\int_{B_r(x)}|u_m|^2\,dx=0.
\]
Then Lemma~\ref{Lem-3.3} implies that for every $p\in(2,\infty)$,
\[
\|u_m\|_{L^p(\mathcal G)}\to0.
\]
Fix such a $p$. By $(F_2)$--$(F_3)$, for every $\varepsilon>0$ there exists $C_\varepsilon>0$ such that
\[
|F_u(x,u)|\le \varepsilon |u| + C_\varepsilon |u|^{p-1}
\qquad\text{for all }(x,u)\in\mathcal G\times\mathbb C^2.
\]
Testing $\Phi'(u_m)$ by $u_m^+-u_m^-$ we have
\[
\Phi'(u_m)[u_m^+-u_m^-]
=\|u_m\|_\omega^2-\int_{\mathcal G}\langle F_u(x,u_m),\,u_m^+-u_m^-\rangle\,dx,
\]
hence, since $\Phi'(u_m)\to0$ in $Y^*$ and $(u_m)$ is bounded in $Y$,
\[
\|u_m\|_\omega^2
=\int_{\mathcal G}\langle F_u(x,u_m),\,u_m^+-u_m^-\rangle\,dx + o(1).
\]
Using the above growth estimate, H\"older's inequality, and $\|u_m\|_{L^p}\to0$, we obtain
\[
\Bigl|\int_{\mathcal G}\langle F_u(x,u_m),\,u_m^+-u_m^-\rangle\,dx\Bigr|
\le \varepsilon C_1 + C_\varepsilon C_2\,\|u_m\|_{L^p(\mathcal G)}^{p-1}\to \varepsilon C_1,
\]
for some constants $C_1,C_2>0$ independent of $m$. Thus
\[
\limsup_{m\to\infty}\|u_m\|_\omega^2 \le \varepsilon C_1.
\]
Since $\varepsilon>0$ is arbitrary, we get $\|u_m\|_\omega\to0$, hence $\|u_m\|\to0$. In particular,
\[
\Phi(u_m)-\frac12\Phi'(u_m)[u_m]
=\int_{\mathcal G}\hat F(x,u_m)\,dx \to 0,
\]
and since $\Phi'(u_m)[u_m]\to0$, we infer $\Phi(u_m)\to0$, i.e. $c=0$.
This contradicts $c\ge\varsigma>0$. Therefore $(u_m)$ is nonvanishing.

\textbf{Step 2.}
By Lemma~\ref{Lem-3.3} and the $\mathbb Z^d$--invariance of $\Phi$, there exist a sequence
$(k_m)\subset\mathbb Z^d$ and a function $v\in Y\setminus\{0\}$ such that, up to a subsequence,
\[
v_m=k_m*u_m\rightharpoonup v \ \text{ in }Y,\qquad
v_m\to v \ \text{ in }L^p_{\mathrm{loc}}(\mathcal G,\mathbb C^2)\ \text{for all }p\in[2,\infty).
\]
Since each translation $k*u$ acts isometrically on $Y$ and $\Phi(k*u)=\Phi(u)$,
we also have the covariance of the derivative:
\[
\Phi'(k*u)[\varphi]=\Phi'(u)[k^{-1}*\varphi]\qquad(\forall\,u,\varphi\in Y,\ k\in\mathbb Z^d),
\]
hence $(v_m)$ is again a $(C)_c$--sequence and $\Phi'(v_m)\to0$ in $Y^*$.

By Lemma~\ref{lem-4.1}, $\Phi'$ is weakly sequentially continuous on $Y$; therefore,
from $v_m\rightharpoonup v$ in $Y$ we obtain
\[
\Phi'(v_m)\to \Phi'(v)\quad\text{in }Y^{*}.
\]
Combining this with $\Phi'(v_m)\to0$ gives $\Phi'(v)=0$, and $v\neq0$ by construction.
By Proposition~\ref{pro-3.1}, $v$ is a bound state of frequency $\omega$ for
NLDE \eqref{eq-1.4}. This proves the existence part of Theorem~\ref{theo-1.1}.

\textbf{Step 3.}
Assume now by contradiction that NLDE \eqref{eq-1.4} admits only finitely many
geometrically distinct bound states. This is equivalent to \eqref{eq-4.22},
with $\mathbb{Z}^{d}$ acting by translations as above.

Under this assumption, Lemma~\ref{lem-4.1}, Lemma~\ref{lem-4.2},
Lemma~\ref{lem-4.3} and Lemma~\ref{lem-4.7} show that $\Phi$ satisfies all
the conditions $(\Phi_{0})$--$(\Phi_{5})$. Moreover, $\Phi$ is even and
$\Phi(0)=0$ by the assumptions on $F$. Hence we may apply
Theorem~\ref{theo-3.2}, which yields an unbounded sequence of positive
critical values of $\Phi$.

This contradicts \eqref{eq-4.22}. Indeed, since $\Phi(k*u)=\Phi(u)$ for all
$k\in\mathbb Z^d$, each $\mathbb Z^d$--orbit of critical points contributes
only one critical value, so finitely many $\mathbb Z^d$--orbits can generate
only a finite set of critical values. Therefore \eqref{eq-4.22} cannot hold,
and NLDE \eqref{eq-1.4} possesses infinitely many geometrically distinct
bound states.

\section*{Acknowledgments}
We would like to thank the anonymous referee for his/her careful readings of our manuscript and the useful comments. 

\medskip
{\bf Funding:} This work is supported by National Natural Science Foundation of China (12301145, 12261107, 12561020) and Yunnan Fundamental Research Projects (202301AU070144, 202401AU070123).

\medskip
{\bf Author Contributions:} All the authors wrote the main manuscript text together and these authors contributed equally to this work.

\medskip
{\bf Data availability:}  Data sharing is not applicable to this article as no new data were created or analyzed in this study.

\medskip
{\bf Conflict of Interests:} The author declares that there is no conflict of interest.


\bibliographystyle{plain}
\bibliography{reference}
\end{document}